\documentclass[12pt]{amsart}
\usepackage{amsmath}
\usepackage{amssymb}
\usepackage[shortlabels]{enumitem}
\usepackage{booktabs}
\usepackage{url}
\usepackage{hyperref}
\usepackage{tikz}
\usepackage{tikz-cd}
\hypersetup{unicode=true,
colorlinks=true,
citecolor=black,
linkcolor=black,
urlcolor=black,
plainpages=false}
\allowdisplaybreaks
\newtheorem{teo}{Theorem}[section]
\newtheorem{theorem}[teo]{Theorem}
\newtheorem{proposition}[teo]{Proposition}
\newtheorem{lemma}[teo]{Lemma}
\newtheorem{corollary}[teo]{Corollary}
\theoremstyle{definition}
\newtheorem{definition}[teo]{Definition}
\newtheorem{example}[teo]{Example}
\theoremstyle{remark}
\newtheorem{remark}[teo]{Remark}

\newcommand{\CC}{\mathcal{C}}
\newcommand{\DD}{\mathcal{D}}
\newcommand{\EE}{\mathcal{E}}
\newcommand{\PP}{\mathcal{P}}
\newcommand{\co}{\colon}
\newcommand{\id}{\mathrm{id}}
\newcommand{\N}{\mathrm{N}}
\newcommand{\Hom}{\mathrm{Hom}}
\newcommand{\Ob}{\mathrm{Ob}}
\newcommand{\Z}{\mathbb{Z}}
\newcommand{\K}{\mathcal{K}}
\newcommand{\Q}{\mathcal{Q}}
\newcommand{\M}{M}
\newcommand{\RP}{\mathbb{RP}}
\newcommand{\cone}{\mathrm{Cone}}
\newcommand{\op}{\mathrm{op}}
\newcommand{\nerve}{\Delta}
\newcommand{\lk}{\mathrm{lk}}
\newcommand{\st}{\mathrm{star}}
\newcommand{\dl}{\mathrm{dl}}
\newcommand{\Sph}{\mathcal{S}}
\newcommand{\tSph}{\widetilde{\mathcal{S}}}

\begin{document}

\title[Minimal acyclic categorical models]{How many points does it take to see the homotopy of a space?\\ Minimal acyclic categorical models}

\author[Isaac Carcac\'ia Campos]{%
Isaac Carcac\'ia Campos}
\address{%
Isaac Carcac\'ia Campos
\\
Departamento de Matem\'aticas, Universidade de Santiago de Compostela, 15782-SPAIN}
\email{isaac.c.campos@usc.es}

\begin{abstract}

Every connected homotopy type is the classifying space of a small category, but this general result gives little control over the categories involved. We study finite acyclic categorical models and measure their complexity by the number of objects.

Our main tool is the sink lemma, which describes the addition of a sink as the attachment of a cone along an explicit map. Using it, we determine the minimum number of objects required to model spheres, bouquets of spheres, Moore spaces of cyclic groups, real projective spaces, and lens spaces. We also exhibit infinitely many pairwise non-isomorphic rigid minimal models of the sphere.

To include spaces admitting no finite model, we introduce a truncated invariant measuring the minimum number of objects required for a finite acyclic category to approximate a space through degree $n$. Spreading a finite monoid over $n+1$ objects produces such a category with an $n$-connected map to the monoid's classifying space. For every nontrivial finite group, this bound is optimal, although its classifying space admits no finite acyclic model.
\end{abstract}


\maketitle

\section*{Introduction}

No branch of mathematics is an island. When the bell tolls for a discovery in one of them, it tolls for all: the map of mathematics is connected, and however far apart two regions may seem, there are paths between them along which a change in one place is heard in the others. Along these paths travel amphibious objects, which live on both shores at once and carry problems from one area to another. Category theory provides a language for such passages. In this paper we use it to move between topological spaces and finite combinatorial structures, with the hope of shedding some light on both shores.

The question behind this paper is a simple one: how many points does it take to represent a space? Finite posets are the first amphibious objects in this story. They are topological spaces, since finite posets are the same as finite $T_0$-spaces \cite{Alexandroff1937}, and they are categories, namely those with at most one morphism between any two objects. McCord proved that a finite $T_0$-space and the classifying space of its poset have the same weak homotopy type \cite{McCord1966}, so the minimum number of points of a finite model becomes a measure of how complicated a space is, and it has been studied in depth \cite{BARMAK,BarmakMinian2007,CianciOttina2018}.

In a poset two points are either related or not. We allow them to be related in several ways, and work with \emph{acyclic categories}: small categories with no morphisms in both directions between distinct objects and no endomorphisms other than identities. Kozlov developed their combinatorial topology systematically \cite[Chapters~10 and~14]{Kozlov2008}, and Tanaka studied their strong homotopy theory \cite{Tanaka2019}. A finite acyclic category $\CC$ models a space through its classifying space $B\CC$, a finite CW complex, and parallel morphisms make a big difference: the Kronecker category, with two objects and two parallel morphisms, has a circle as classifying space, whereas a poset needs four points to model $S^1$.

Acyclicity is what keeps the count honest. Without it a single object is enough: every connected CW complex is homotopy equivalent to the classifying space of a discrete monoid \cite{McDuff1979}, and by Ben-Bassat's recent proof of a conjecture of Fiedorowicz \cite[p.~323]{Fiedorowicz1984}, every finite simply connected CW complex is, up to homotopy, the classifying space of a finite monoid \cite{BenBassat2026}. In a monoid all the complexity is hidden in the endomorphisms of a single point. In an acyclic category there are no endomorphisms to hide it in, and the complexity has to be shared between the objects and the parallel morphisms. We measure it with two invariants. The \emph{object complexity} $\nu(X)$ of a space $X$ is the smallest number of objects of a finite acyclic category $\CC$ with $B\CC\simeq X$, and the \emph{truncated object complexity} $\nu_n(X)$ is the smallest number of objects of one that models $X$ through degree $n$, that is, that admits an $n$-connected map $B\CC\to X$. The numbers $\nu_n(X)$ increase with $n$ and are bounded by $\nu(X)$. For a finite complex they eventually reach $\nu(X)$; for a space with no finite model they measure the price of seeing it one degree at a time.

Almost all of our models are understood through a single construction. A CW complex is built by attaching cells one at a time, and a finite acyclic category can be built in the same way, adding one object at a time, each new object being a sink of what has been built so far. Adding a sink glues a cone to the classifying space along a map that can be read off from the morphisms into the sink: this is the \emph{sink lemma} (Lemma~\ref{lem:sink}). Up to homotopy it is a special case of a recent collage lemma of Ben-Bassat \cite[Lemma~2.5]{BenBassat2026}; in the acyclic setting it holds at the level of nerves, and this is what allows us to compute attaching maps and their degrees. With three objects, it reduces the homotopy type of a category to a map from a graph to a wedge of circles, which is determined by homology when there are two morphisms between the first two objects (Proposition~\ref{prop:three}).

Both invariants obey simple homological lower bounds: a nonzero class in degree $k$ costs $k+1$ objects, and torsion costs one more (Proposition~\ref{prop:bounds}). These bounds are attained by spheres and bouquets of spheres of the same dimension, by the Moore spaces of cyclic groups, and by real projective spaces and lens spaces, for which we construct explicit object-minimal models. Classifying spaces of finite groups have no finite model, but their truncated invariants are as small as the bounds allow. Spreading a finite monoid $M$ over $n+1$ objects gives a category whose classifying space maps $n$-connectedly to $BM$ (Theorem~\ref{thm:monoids}), and for every nontrivial finite group $G$ this gives $\nu_n(BG)=n+1$ for every $n$ (Proposition~\ref{prop:groups-minimal}); already with three objects one obtains every finite group as a fundamental group. Combined with Ben-Bassat's theorem, the same construction gives $\nu_n(X)\le n+1$ for every finite simply connected complex $X$.

The last difference with posets concerns uniqueness. The minimal finite model of $S^n$ is unique \cite[Corollary~2.14]{BarmakMinian2007}, but $S^2$ has infinitely many object-minimal models with three objects (Corollary~\ref{cor:infinite}). They have no beat objects in the sense of Tanaka, so they are rigid: two of them are strongly homotopy equivalent only if they are isomorphic \cite{Tanaka2019}. Even among the most rigid categories, the homotopy type of the classifying space does not determine the category. A space can be seen through few points in many essentially different ways.

The paper is organized as follows. Section~\ref{sec:prelim} collects background on acyclic categories, their nerves and beat objects. Section~\ref{sec:sink} proves the sink lemma and studies categories with three objects. Section~\ref{sec:bounds} introduces $\nu$ and $\nu_n$ together, with their lower bounds, the stabilization of $\nu_n$ at $\nu$ for finite complexes, and the relation between object-minimality and strong minimality. Sections~\ref{sec:bouquets}, \ref{sec:moore} and~\ref{sec:projective} compute $\nu$ for spheres and bouquets, Moore spaces, and real projective and lens spaces, and Section~\ref{sec:monoids} turns to the truncated invariant and to classifying spaces of monoids and groups. Section~\ref{sec:comparison} compares our models with minimal finite models and returns to the question of where the complexity of a space can be stored.

\subsection*{Declaration on the use of generative AI}
Claude Opus 5.5 was used to assist in preparing diagrams, improving the English exposition, translating passages from Spanish into English, shortening the manuscript, and identifying potential gaps or errors in the proofs and calculations. All AI-assisted material was reviewed and independently verified by the author, who assumes full responsibility for the mathematical content. The AI also suggested removing several digressions and philosophical reflections; those that remain do so at the author’s insistence.
\section{Preliminaries}\label{sec:prelim}

This section collects the background, all of it known, on which the rest of the paper is built. We follow Kozlov's treatment of acyclic categories and trisps \cite[Section~2.3 and Chapter~10]{Kozlov2008} and Tanaka's theory of beat objects \cite{Tanaka2019}.

\subsection{Acyclic categories and their nerves}\label{subsec:acyclic}\label{subsec:nerve}

A small category is \emph{acyclic} if its only isomorphisms and its only endomorphisms are identities; equivalently, for distinct objects $x,y$ at most one of $\Hom(x,y)$ and $\Hom(y,x)$ is nonempty, and $\Hom(x,x)=\{\id_x\}$ \cite[Section~10.1.1]{Kozlov2008}. A morphism is \emph{non-trivial} if it is not an identity. The relation $x\preceq y$ if $\Hom(x,y)\neq\emptyset$ is a partial order on $\Ob\CC$, which forgets how many morphisms there are from $x$ to $y$ and extends to a total order \cite[Section~10.1.2]{Kozlov2008}; posets are exactly the acyclic categories with at most one morphism between any two objects. A finite nonempty acyclic category has a \emph{source}, an object with no non-trivial morphism into it, and a \emph{sink}, an object with no non-trivial morphism out of it.

For an acyclic category the degenerate simplices of the nerve can be discarded altogether. Following Kozlov, we use \emph{trisps}, also known as semisimplicial sets or $\Delta$-complexes \cite[Section~2.3.1, Definition~2.47]{Kozlov2008}. The \emph{nerve} $\nerve(\CC)$ of an acyclic category $\CC$ is the regular trisp whose vertices are the objects of $\CC$ and whose $k$-simplices are the strings $x_0\xrightarrow{f_1}x_1\to\cdots\xrightarrow{f_k}x_k$ of composable non-trivial morphisms; the face opposite $x_i$ is obtained by composing $f_{i+1}\circ f_i$ if $0<i<k$, and by deleting $f_1$ or $f_k$ if $i=0$ or $i=k$ \cite[Definition~10.4]{Kozlov2008}. Its simplices are the nondegenerate simplices of the simplicial nerve $\N\CC$, with the same faces, so the \emph{classifying space} $B\CC=|\nerve(\CC)|=|\N\CC|$ is a CW complex with one cell for each simplex of $\nerve(\CC)$ \cite[Section~10.2.1]{Kozlov2008}. Functors induce trisp maps \cite[Proposition~10.6]{Kozlov2008}, and natural transformations induce homotopies of classifying spaces \cite{Quillen1973}.

\begin{lemma}\label{lem:faces}
Let $\CC$ be a finite acyclic category.
\begin{enumerate}[(a)]
\item $\dim B\CC\le|\Ob\CC|-1$.
\item The $k+1$ faces of a $k$-simplex of $\nerve(\CC)$ are pairwise distinct, so each column of the cellular boundary $\partial_k$ has exactly $k+1$ nonzero entries, all equal to $\pm1$. For instance, a $2$-simplex $x\xrightarrow{f}y\xrightarrow{g}z$ has boundary $g-g\circ f+f$.
\item If $d=\dim B\CC$, then $H_d(B\CC;\Z)$ is free abelian.
\end{enumerate}
\end{lemma}

\begin{proof}
The objects of a string of non-trivial morphisms are strictly increasing for $\preceq$, which gives (a). For (b), the faces of a simplex have different sets of objects. For (c), there are no $(d+1)$-cells, so $H_d$ is the group of cellular $d$-cycles, a subgroup of a free abelian group.
\end{proof}

We use the edge-path presentation of the fundamental group of a CW complex whose $2$-cells are attached along edge paths \cite[Propositions~1.26 and~1A.2]{Hatcher2002}. For $B\CC$, choose a maximal tree $T$ in the $1$-skeleton; then $\pi_1(B\CC)$ is generated by the classes $[e]$ of the edges, with relations $[e]=1$ for $e\in T$ and $[f][g]=[g\circ f]$ for each $2$-simplex $(f,g)$.

\subsection{Joins, products and source wedges}\label{subsec:constructions}\label{subsec:wedge}

The \emph{join} $\CC\star\DD$ of acyclic categories, called the stack $\CC\oplus\DD$ in \cite[Definition~10.7]{Kozlov2008}, has objects $\Ob\CC\sqcup\Ob\DD$, the morphisms of $\CC$ and of $\DD$, and exactly one morphism $x\to y$ for each $x\in\Ob\CC$ and $y\in\Ob\DD$. We write $\CC^{\star k}$ for the $k$-fold join, $\mathbf n$ for the discrete category with $n$ objects, $\EE^{\triangleright}=\EE\star\mathbf1$, and call $\mathbf 2\star\CC$ the \emph{categorical suspension} of $\CC$; the dual construction $\CC\star\mathbf2$ is the categorical suspension of \cite{CarcaciaCampos2026OpenCovers}. If $\CC$ and $\DD$ are finite with sources $s$ and $s'$, their \emph{source wedge} $\CC\vee\DD$ is obtained from the coproduct $\CC\sqcup\DD$ \cite[Section~10.3.1]{Kozlov2008} by identifying $s$ with $s'$.

\begin{lemma}\label{lem:join}
Let $\CC$ and $\DD$ be acyclic categories.
\begin{enumerate}[(a)]
\item $\nerve(\CC\star\DD)\cong\nerve(\CC)*\nerve(\DD)$, so $B(\CC\star\DD)\cong B\CC*B\DD$ \cite[(10.1), Section~2.3.2]{Kozlov2008}. In particular $B(\EE^{\triangleright})\cong\cone(B\EE)$, $B(\mathbf2\star\CC)\cong\Sigma B\CC$, and the stack of $n$ copies of $\mathbf2$ has classifying space $S^{n-1}$.
\item The product $\CC\times\DD$ is acyclic, and $B(\CC\times\DD)\cong B\CC\times B\DD$ \cite[Section~10.5.1, Theorem~10.21]{Kozlov2008}.
\item If $\CC$ and $\DD$ are finite, the source wedge $\CC\vee\DD$ is a finite acyclic category with $|\Ob\CC|+|\Ob\DD|-1$ objects, and $B(\CC\vee\DD)\cong B\CC\vee B\DD$. The same holds with sinks in place of sources.
\end{enumerate}
\end{lemma}

\begin{proof}
For (c), no non-trivial morphism enters $s$, and there are no morphisms between $\CC\setminus s$ and $\DD\setminus s$, so every string of non-trivial morphisms lies in $\CC$ or in $\DD$, and $\nerve(\CC\vee\DD)$ is the union of $\nerve(\CC)$ and $\nerve(\DD)$ along the vertex $s$. For sinks, pass to opposite categories.
\end{proof}

We use freely the following standard facts for finite CW complexes \cite{Hatcher2002}: the join is homotopy invariant and associative, $X*S^q\cong\Sigma^{q+1}X$, and the unreduced suspension commutes with wedges up to homotopy, so that $\Sigma^j\bigl(\bigvee^rS^p\bigr)\simeq\bigvee^rS^{p+j}$ and $\Sigma^j(S^p\cup_fe^{p+1})\simeq S^{p+j}\cup_{\Sigma^jf}e^{p+j+1}$.

\subsection{Beat objects and strong minimality}\label{subsec:beat}

\begin{definition}[{\cite[Definitions~2.3 and~2.10]{Tanaka2019}}]\label{def:beat}
An object $x$ of an acyclic category $\CC$ is a \emph{down beat object} if there is a non-trivial morphism $f\co y\to x$ such that every non-trivial morphism $g\co z\to x$ factors uniquely as $g=f\circ\tilde g$ with $\tilde g\co z\to y$, possibly an identity. Up beat objects are defined dually. Tanaka calls $\CC$ \emph{minimal} if it has no beat objects; to avoid confusion with the object-minimal models of Section~\ref{sec:bounds}, we call such categories \emph{strongly minimal}.
\end{definition}

For finite posets, beat objects are Stong's beat points \cite{Stong1966}; see \cite[Remark~2.4]{Tanaka2019}, and see \cite[Section~6]{CarcaciaCampos2026Bundles} for relative versions. Two functors are \emph{homotopic} if they are connected by a finite zigzag of natural transformations, and a \emph{strong homotopy equivalence} is a functor with a homotopy inverse in this sense \cite[Definition~2.2]{Tanaka2019}; it induces a homotopy equivalence of classifying spaces.

\begin{lemma}[{\cite[Corollary~2.7, Proposition~2.12]{Tanaka2019}}]\label{lem:tanaka}
Let $\CC$ be an acyclic category.
\begin{enumerate}[(a)]
\item If $x$ is a beat object of $\CC$, then $\CC$ and $\CC\setminus x$ are strongly homotopy equivalent. In particular $B\CC\simeq B(\CC\setminus x)$.
\item If $\CC$ is strongly minimal and $F\co\CC\to\CC$ is homotopic to $1_\CC$, then $F=1_\CC$. Consequently, strongly minimal acyclic categories that are strongly homotopy equivalent are isomorphic.
\end{enumerate}
\end{lemma}

\section{Sinks and mapping cones}\label{sec:sink}

Finite acyclic categories can be built one object at a time, as CW complexes are built one cell at a time. In this section we show that each new object behaves like a cell: adding a sink to a category glues a cone to its classifying space, along a map that is read off from the morphisms into the new object. We then characterize the sinks that are beat objects, and we end with categories with three objects, in which everything is encoded in a map between two graphs.

\subsection{Adding a sink}\label{subsec:add-sink}

Let $\DD$ be a small category and let $P\co\DD^{\op}\to\mathbf{Set}$ be a presheaf. \emph{Adding a sink to $\DD$ along $P$} produces the category $\DD[P]$ with objects $\Ob\DD\sqcup\{t\}$, the morphisms of $\DD$, the elements of $P(x)$ as morphisms $x\to t$, only the identity at $t$, and no morphisms out of $t$. The composite of $f\co x\to y$ in $\DD$ and $u\in P(y)$ is $u\circ f=P(f)(u)$, and associativity is the functoriality of $P$. The presheaf records how the new object sees the old ones: which morphisms reach it from each object, and how these are related by the morphisms of $\DD$. If $\DD$ is finite and acyclic and every $P(x)$ is finite, then $\DD[P]$ is finite and acyclic, and $t$ is a sink. Conversely, if $t$ is a sink of a small category $\CC$ and $\DD=\CC\setminus t$, then $\CC\cong\DD[\Hom_\CC(-,t)|_\DD]$. Hence every finite acyclic category grows in this way from a single object: if $x_0,\dots,x_n$ is a total order on $\Ob\CC$ extending $\preceq$, then $x_k$ is a sink of the full subcategory $\CC_k$ on $\{x_0,\dots,x_k\}$, and $\CC_k\cong\CC_{k-1}[\Hom_\CC(-,x_k)|_{\CC_{k-1}}]$.

From now on let $\CC$ be a finite acyclic category, let $t$ be a sink, let $\DD=\CC\setminus t$ and let $P=\Hom_\CC(-,t)|_\DD$. The \emph{category of elements} of $P$ is the category $\DD/t$ whose objects are the pairs $(x,u)$ with $x\in\Ob\DD$ and $u\co x\to t$, and whose morphisms $(x,u)\to(y,v)$ are the morphisms $f\co x\to y$ of $\DD$ with $v\circ f=u$. It is the category $\CC_{>t}$ of \cite[Section~10.3.3]{Kozlov2008} and Tanaka's category $\widehat{\CC/t}$ \cite[Section~2]{Tanaka2019}. It is finite and acyclic, and we write $\mathrm{pr}\co\DD/t\to\DD$ for the projection $(x,u)\mapsto x$.

For a vertex $x$ of $\nerve(\CC)$, let $\dl(x)$, $\lk(x)$ and $\st(x)=\cone(\lk(x))$ be its deletion, link and star \cite[Section~2.3.2]{Kozlov2008}. Since there are no non-trivial morphisms out of $t$, Kozlov's description of links \cite[Section~10.3.3, (10.2), (10.3)]{Kozlov2008} gives $\lk(t)\cong\nerve(\DD/t)$ and $\st(t)\cong\cone(\nerve(\DD/t))$, while $\dl(t)=\nerve(\DD)$. The star is not embedded in $\nerve(\CC)$ in general: on the base of the cone, the map to $\nerve(\CC)$ is $\nerve(\mathrm{pr})$, which may identify simplices. The sink lemma says that this is the only identification.

\begin{lemma}[Sink lemma]\label{lem:sink}
In the situation above:
\begin{enumerate}[(a)]
\item The trisp $\nerve(\CC)$ is obtained from $\nerve(\DD)$ by attaching the cone $\nerve((\DD/t)^{\triangleright})=\cone(\nerve(\DD/t))$ along the trisp map $\nerve(\mathrm{pr})\co\nerve(\DD/t)\to\nerve(\DD)$. Consequently $B\CC$ is homeomorphic to the mapping cone
\[
B\DD\cup_{B\mathrm{pr}}\cone B(\DD/t).
\]
If $\DD/t$ is empty, then $t$ is an isolated object and $B\CC=B\DD\sqcup\{t\}$.
\item The inclusion $B(\DD/t)\to\cone B(\DD/t)$ is the inclusion of a subcomplex, so the pushout in (a) is a homotopy pushout of $B\DD\leftarrow B(\DD/t)\to\ast$.
\item If $\DD/t\neq\emptyset$, there is a long exact sequence
\begin{multline*}
\cdots\to\widetilde H_k(B(\DD/t))\xrightarrow{\mathrm{pr}_*}\widetilde H_k(B\DD)\to\widetilde H_k(B\CC)\\
\to\widetilde H_{k-1}(B(\DD/t))\xrightarrow{\mathrm{pr}_*}\cdots.
\end{multline*}
\item If $B\DD$ and $B(\DD/t)$ are nonempty and connected, choose an object $(x_0,u_0)$ of $\DD/t$. Then
\[
\pi_1(B\CC,x_0)\cong\pi_1(B\DD,x_0)\big/\bigl\langle\!\bigl\langle(B\mathrm{pr})_*\pi_1\bigl(B(\DD/t),(x_0,u_0)\bigr)\bigr\rangle\!\bigr\rangle.
\]
\end{enumerate}
\end{lemma}

\begin{proof}
(a) Let $J\co(\DD/t)^{\triangleright}\to\CC$ be the functor with $J(x,u)=x$, $J(f)=f$, $J(*)=t$ for the adjoined terminal object $*$, and sending the unique morphism $(x,u)\to*$ to $u$. It is a functor: for $f\co(x,u)\to(y,v)$, the composite $(x,u)\to(y,v)\to*$ is the unique morphism $(x,u)\to*$, and $J$ sends both to $v\circ f=u$. Its restriction to $\DD/t$ is $\mathrm{pr}$ followed by the inclusion $\DD\subseteq\CC$. By Lemma~\ref{lem:join}, $\nerve((\DD/t)^{\triangleright})=\cone(\nerve(\DD/t))$, whose simplices are those of $\nerve(\DD/t)$, the apex, and the cones over simplices of $\nerve(\DD/t)$.

The simplices of $\nerve(\CC)$ not containing $t$ are those of $\nerve(\DD)$, since $\DD$ is full. Since $t$ is a sink, a simplex containing $t$ has $t$ as its last vertex, and is either the vertex $t$ or a string
\[
x_0\xrightarrow{f_1}\cdots\xrightarrow{f_{p-1}}x_{p-1}\xrightarrow{u}t
\]
with $x_0,\dots,x_{p-1}\in\Ob\DD$. Such a string is the image under $\nerve(J)$ of exactly one simplex of the cone not in its base, namely the cone over
\[
(x_0,u_0)\xrightarrow{f_1}\cdots\xrightarrow{f_{p-1}}(x_{p-1},u_{p-1}),\qquad u_i=u\circ f_{p-1}\circ\cdots\circ f_{i+1},
\]
which is a string in $\DD/t$ because $u_{i+1}\circ f_{i+1}=u_i$; any preimage has these underlying morphisms, and its objects are forced by $u_{p-1}=u$ and these compatibilities. The vertex $t$ is the image of the apex. Thus $\nerve(J)$ and the inclusion of $\nerve(\DD)$ give a bijection between the simplices of $\nerve(\DD)\sqcup\bigl(\cone(\nerve(\DD/t))\setminus\nerve(\DD/t)\bigr)$ and those of $\nerve(\CC)$, compatible with faces: deleting a vertex other than $t$ gives the corresponding face in the cone, and deleting $t$ gives the image under $\nerve(\mathrm{pr})$ of the base simplex. This is the gluing data of the pushout of $\nerve(\DD)\xleftarrow{\nerve(\mathrm{pr})}\nerve(\DD/t)\hookrightarrow\cone(\nerve(\DD/t))$. Realization of trisps is obtained by gluing simplices along faces, so it preserves this pushout, and $|\cone(\nerve(\DD/t))|$ is the cone on $B(\DD/t)$ by Lemma~\ref{lem:join}. If $\DD/t$ is empty, no non-trivial morphism has target $t$, and the cone is the single vertex $t$.

(b) The inclusion of a subcomplex is a cofibration, and the cone is contractible. (c) If $Z$ is the mapping cylinder of $B\mathrm{pr}$, then $Z\simeq B\DD$ and $Z/B(\DD/t)$ is the mapping cone, so (c) is the long exact sequence of the pair $(Z,B(\DD/t))$. (d) Apply the van Kampen theorem to the open sets of the mapping cone given by the images of $B\DD\sqcup B(\DD/t)\times[0,2/3)$, which deformation retracts onto $B\DD$, and of $B(\DD/t)\times(1/3,1]$, which is contractible; their intersection deformation retracts onto $B(\DD/t)$, and after moving the base point along $\{(x_0,u_0)\}\times[0,1/2]$ its inclusion into the first set induces $(B\mathrm{pr})_*$.
\end{proof}

The first examples are the familiar ones. If $P$ is the terminal presheaf, then $\DD[P]=\DD^{\triangleright}$, $\DD/t\cong\DD$ and $\mathrm{pr}$ is the identity, so the sink lemma gives the cone on $B\DD$. Adding to $\DD^{\triangleright}$ a second sink along the presheaf that is a point on $\DD$ and empty at the apex gives $\DD\star\mathbf 2$, and the suspension of $B\DD$ appears as the union of two cones. The interesting cases are those in which $\mathrm{pr}$ is far from injective. Then the attaching map can fold several pieces onto one, or wrap around a circle several times, and this is where the economy of acyclic categories comes from.

\begin{remark}\label{rem:collage}
Up to homotopy, the sink lemma is a special case of the collage lemma of Ben-Bassat: the collage of a functor $Q\co\mathcal A^{\op}\times\mathcal B\to\mathbf{Set}$ has classifying space homotopy equivalent to the homotopy pushout of $B\mathcal A\leftarrow B\,\mathrm{El}(Q)\to B\mathcal B$ \cite[Definition~2.3, Lemma~2.5]{BenBassat2026}, and for $\mathcal B=\mathbf 1$ the collage is $\DD[P]$ and $\mathrm{El}(Q)=\DD/t$. What the acyclic setting adds is a statement at the level of nerves: $\nerve(\CC)$ is obtained from $\nerve(\DD)$ by attaching a cone along the explicit trisp map $\nerve(\mathrm{pr})$, so that $B\CC$ is homeomorphic, and not only homotopy equivalent, to the mapping cone, and the attaching map can be read off simplex by simplex. This is what we use to compute the degrees of attaching maps in Sections~\ref{sec:bouquets} and~\ref{sec:moore}. Both statements carry to small categories the splitting technique of Cianci and Ottina, who cut a finite poset into a subspace and its complement to control its homotopy type and the size of its finite models \cite[Section~3]{CianciOttina2018}.
\end{remark}

\subsection{Representable presheaves and beat objects}\label{subsec:representable}

A sink need not be a terminal object: the sets $P(x)$ may be empty or have several elements. The simplest sinks are those added along representable presheaves, and they are exactly the sinks that are beat objects.

\begin{lemma}\label{lem:beat-representable}
Let $\CC$ be an acyclic category, let $x$ be an object of $\CC$, and let $\CC_{<x}$ be the full subcategory on the objects $y\neq x$ with $\Hom_\CC(y,x)\neq\emptyset$. Then $x$ is a down beat object if and only if the presheaf $\Hom_\CC(-,x)|_{\CC_{<x}}$ is representable. Dually, $x$ is an up beat object if and only if $\Hom_\CC(x,-)$, restricted to the full subcategory $\CC_{>x}$ of the objects $y\neq x$ with $\Hom_\CC(x,y)\neq\emptyset$, is corepresentable. In particular, a sink $t$ of $\CC$ is a beat object if and only if the presheaf $\Hom_\CC(-,t)|_\DD$ on $\DD=\CC\setminus t$ is representable.
\end{lemma}

\begin{proof}
Let $Q=\Hom_\CC(-,x)|_{\CC_{<x}}$. By the Yoneda lemma, an element $f\co y\to x$ of $Q(y)$ determines the natural transformation $\Hom_{\CC_{<x}}(-,y)\to Q$ given by $\tilde g\mapsto f\circ\tilde g$, and every natural transformation between these presheaves is of this form. Every non-trivial morphism into $x$ has its source in $\CC_{<x}$, and $\CC_{<x}$ is full, so the condition of Definition~\ref{def:beat} says precisely that this transformation is bijective at every object of $\CC_{<x}$: surjectivity is the existence of the factorization $g=f\circ\tilde g$, and injectivity is its uniqueness. The dual statement follows by passing to opposite categories. For the last assertion, a sink can only be a down beat object, and $P=\Hom_\CC(-,t)|_\DD$ is representable on $\DD$ if and only if its restriction to $\CC_{<t}$ is: if $f\in P(y)$ and $\DD(z,y)\neq\emptyset$, then $P(z)\neq\emptyset$, so both $P$ and $\DD(-,y)$ vanish outside $\CC_{<t}$.
\end{proof}

The sink lemma shows directly why such sinks do not change the homotopy type. If $P$ is represented by $y$, with universal element $f\co y\to t$, then $(y,f)$ is a terminal object of $\DD/t$: for every object $(z,u)$ there is a unique $\tilde g\co z\to y$ with $f\circ\tilde g=u$, and such a $\tilde g$ is precisely a morphism $(z,u)\to(y,f)$. Hence $B(\DD/t)$ is contractible, and by Lemma~\ref{lem:sink}(b) the inclusion $B\DD\to B\CC$ is a homotopy equivalence; this is the shadow, on classifying spaces, of Lemma~\ref{lem:tanaka}(a). The converse fails, and the gap between the two notions is the gap between being contractible and having a terminal object. For instance, if $\DD$ is the poset with $a<b$ and $a<c$ and $\CC=\DD^{\triangleright}$, then $\DD/t\cong\DD$ has the initial object $(a,u_a)$, so $B(\DD/t)$ is contractible and $B\DD\to B\CC$ is a homotopy equivalence; but $P$ is the terminal presheaf on $\DD$, which is not representable since $\DD$ has no terminal object, so $t$ is not a beat object.

\subsection{Three objects}\label{subsec:three}

Three objects are already enough to see the sink lemma at work, and to see how much information a composition table can hide. Let $\CC$ be an acyclic category with objects $0,1,2$ and no non-trivial morphisms from a larger object to a smaller one, and write
\begin{gather*}
\Hom(0,1)=\{f_0,\dots,f_{a-1}\},\qquad\Hom(1,2)=\{g_0,\dots,g_{b-1}\},\\
\Hom(0,2)=\{h_0,\dots,h_{c-1}\}.
\end{gather*}
Since there are no strings of three non-trivial morphisms, $\CC$ is determined by these sets and by its \emph{composition table}, the $a\times b$ matrix $T$ with $g_j\circ f_i=h_{T(i,j)}$, and every table defines a category. We call $(a,b)$ the \emph{multiplicities} of $\CC$.

\begin{figure}[ht]
\centering
\begin{minipage}[b]{0.45\textwidth}
\centering
\begin{tikzpicture}[scale=1.1]
\coordinate (h) at (0,0);
\coordinate (g0) at (-1.7,0);
\coordinate (g1) at (1.7,0);
\draw[thick] (h) to[bend right=40] node[above]{\scriptsize $e_{00}$} (g0);
\draw[dashed, thick] (h) to[bend left=40] node[below]{\scriptsize $e_{10}$} (g0);
\draw[thick] (h) to[bend left=40] node[above]{\scriptsize $e_{01}$} (g1);
\draw[dashed, thick] (h) to[bend right=40] node[below]{\scriptsize $e_{11}$} (g1);
\foreach \p in {(h),(g0),(g1)} \fill \p circle (1.4pt);
\node[above] at (h) {\small $h$};
\node[left] at (g0) {\small $g_0$};
\node[right] at (g1) {\small $g_1$};
\end{tikzpicture}

\smallskip
(a) $\Gamma_T$ for $\Sph$
\end{minipage}
\hfill
\begin{minipage}[b]{0.45\textwidth}
\centering
\begin{tikzpicture}[scale=1.1]
\coordinate (h0) at (0,0);
\coordinate (g0) at (-1.7,0);
\coordinate (g1) at (1.7,0);
\coordinate (g2) at (0,-1.1);
\coordinate (h1) at (0,-2.1);
\draw[thick] (h0) to[bend right=40] node[above]{\scriptsize $e_{00}$} (g0);
\draw[dashed, thick] (h0) to[bend left=40] node[below]{\scriptsize $e_{10}$} (g0);
\draw[thick] (h0) to[bend left=40] node[above]{\scriptsize $e_{01}$} (g1);
\draw[dashed, thick] (h0) to[bend right=40] node[below]{\scriptsize $e_{11}$} (g1);
\draw[thick] (h0) -- node[right]{\scriptsize $e_{02}$} (g2);
\draw[dashed, thick] (h1) -- node[right]{\scriptsize $e_{12}$} (g2);
\foreach \p in {(h0),(g0),(g1),(g2),(h1)} \fill \p circle (1.4pt);
\node[above] at (h0) {\small $h_0$};
\node[left] at (g0) {\small $g_0$};
\node[right] at (g1) {\small $g_1$};
\node[left] at (g2) {\small $g_2$};
\node[left] at (h1) {\small $h_1$};
\end{tikzpicture}

\smallskip
(b) $\Gamma_T$ for $\tSph_3$
\end{minipage}
\caption{The graphs $\Gamma_T$ of the categories $\Sph$ and $\tSph_3$ of Section~\ref{sec:bouquets}. Solid edges lie over $f_0$ and dashed edges over $f_1$.}
\label{fig:graphs}
\end{figure}
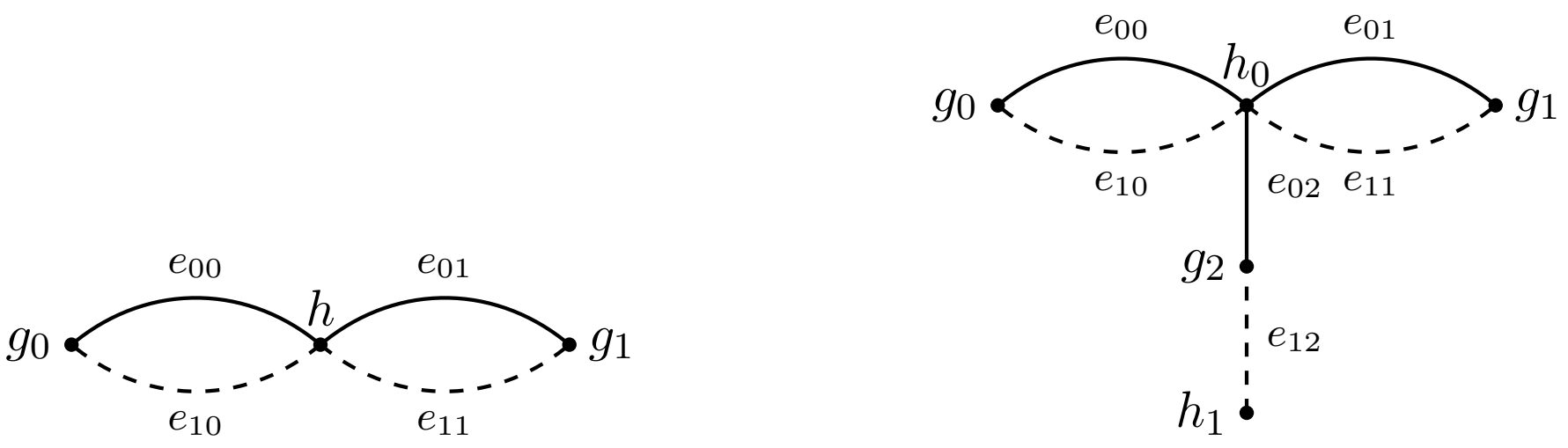

We apply the sink lemma to $t=2$. The full subcategory on $\{0,1\}$ is $\K_a$ (Section~\ref{sec:bouquets}), whose classifying space is a graph with two vertices and $a$ edges, written $B\K$ when $a=2$. The objects of $\DD/2$ are the pairs $(0,h_k)$ and $(1,g_j)$, and its non-trivial morphisms are the $f_i\co(0,h_{T(i,j)})\to(1,g_j)$. There are no composable pairs, so $B(\DD/2)$ is a graph $\Gamma_T$, with an edge $e_{ij}$ from $h_{T(i,j)}$ to $g_j$ that $p_T=B\mathrm{pr}$ sends to the edge $f_i$, and
\begin{equation}\label{eq:three}
B\CC\cong B\K_a\cup_{p_T}\cone(\Gamma_T).
\end{equation}
If $\Gamma_T$ is connected, the rank of $H_1(\Gamma_T;\Z)$ is $ab-b-c+1$. When $a=2$, maps from a graph to the circle $B\K$ are classified by homology, and this gives a complete answer (see Figure~\ref{fig:graphs}).

\begin{proposition}\label{prop:three}
Let $\CC$ be as above with $a=2$, and suppose that $\Gamma_T$ is connected and is not a tree. Let $n\ge1$ be the rank of $H_1(\Gamma_T;\Z)$, and let $d\ge0$ be such that $(p_T)_*\co H_1(\Gamma_T;\Z)\to H_1(B\K;\Z)\cong\Z$ has image $d\Z$. Then
\[
B\CC\simeq\bigl(S^1\cup_de^2\bigr)\vee\bigvee^{n-1}S^2,
\]
which is $\bigvee^{n-1}S^2$, $\M(\Z/d,1)\vee\bigvee^{n-1}S^2$ or $S^1\vee\bigvee^nS^2$ according as $d=1$, $d\ge2$ or $d=0$.
\end{proposition}

\begin{proof}
By \eqref{eq:three} and Lemma~\ref{lem:sink}(b), $B\CC$ is the homotopy pushout of $S^1\leftarrow\Gamma_T\to*$. Collapsing a maximal tree, $\Gamma_T$ is homotopy equivalent to a wedge $W$ of $n$ circles \cite[Proposition~0.17]{Hatcher2002}, so $B\CC$ is homotopy equivalent to the mapping cone of a map $q\co W\to S^1$ \cite[Proposition~0.18 and Section~4.G]{Hatcher2002}. Up to homotopy, $q$ is determined by its degrees $(\delta_1,\dots,\delta_n)$ on the circles, which generate $d\Z$, and precomposing $q$ with self-equivalences of $W$ performs elementary operations on them. By the Euclidean algorithm we may assume that the degrees are $(d,0,\dots,0)$, and then the mapping cone is $(S^1\cup_de^2)\vee\bigvee^{n-1}S^2$.
\end{proof}

\begin{corollary}\label{cor:three-S2}
Let $\CC$ be a category with three objects as above, with $a=2$. Then $B\CC\simeq S^2$ if and only if $\Gamma_T$ is connected, $c=b-1$, and $(p_T)_*\co H_1(\Gamma_T;\Z)\to H_1(B\K;\Z)$ is surjective.
\end{corollary}

\begin{proof}
Suppose that $B\CC\simeq S^2$. Then $\Gamma_T$ is nonempty, since otherwise $2$ is an isolated object, and it is connected, since otherwise the exact sequence of Lemma~\ref{lem:sink}(c) maps $H_1(B\CC;\Z)$ onto $\widetilde H_0(\Gamma_T;\Z)\neq0$. It is not a tree, since otherwise $B\CC\simeq S^1$ by \eqref{eq:three}. By Proposition~\ref{prop:three}, $H_1(B\CC;\Z)=0$ forces $d=1$, and then $n=2b-b-c+1=2$, that is, $c=b-1$. The converse is Proposition~\ref{prop:three} with $n=2$ and $d=1$.
\end{proof}

\section{Object complexity and truncated object complexity}\label{sec:bounds}

We now introduce the two invariants that organize the rest of the paper. The first counts the objects of the smallest exact acyclic model of a space, the second those of the smallest model that is correct through a given degree. Both obey the same homological lower bounds, the second converges to the first for finite complexes, and the rest of the paper shows that the bounds are often attained. Nothing is free, however: what is saved in objects is paid in morphisms.

\subsection{The invariants}\label{subsec:nu}\label{subsec:nu-n}

A map $f\co Y\to X$ is \emph{$n$-connected} if $f_*\co\pi_k(Y,y)\to\pi_k(X,f(y))$ is bijective for $k<n$ and surjective for $k=n$, for every base point $y$ \cite[Section~4.1]{Hatcher2002}. Such a map induces isomorphisms on $H_k$ for $k<n$ and a surjection on $H_n$ \cite[Section~4.2]{Hatcher2002}.

\begin{definition}\label{def:nu}\label{def:nu-n}
The \emph{object complexity} of a space $X$ is
\[
\nu(X)=\min\{|\Ob\CC| : \CC \text{ a finite acyclic category with } B\CC\simeq X\}
\]
if such categories exist, and $\nu(X)=\infty$ otherwise. A finite acyclic category $\CC$ with $B\CC\simeq X$ and $|\Ob\CC|=\nu(X)$ is an \emph{object-minimal model} of $X$. If $X$ is connected and $n\ge0$, a finite acyclic category $\CC$ \emph{models $X$ through degree $n$} if there is an $n$-connected map $B\CC\to X$, and the \emph{truncated object complexity} $\nu_n(X)$ is the smallest number of objects of such a category, with $\nu_n(X)=\infty$ if there is none.
\end{definition}

An $(n+1)$-connected map is $n$-connected, and a homotopy equivalence is $n$-connected for every $n$, so $\nu_n(X)\le\nu_{n+1}(X)\le\nu(X)$ for every connected space $X$ and every $n\ge0$.

\begin{remark}\label{rem:nu-finite}
The invariant $\nu(X)$ is finite exactly when $X$ has the homotopy type of a finite CW complex: one direction is clear, and conversely every finite CW complex is homotopy equivalent to a finite simplicial complex \cite[Theorem~2C.5]{Hatcher2002}, which is modeled by its face poset \cite[(10.4), (10.5)]{Kozlov2008}. Since the barycentric subdivision of a finite acyclic category is a poset with the subdivided nerve \cite[Definition~10.10, (10.7)]{Kozlov2008}, $\nu(X)$ is at most the minimum number of points of a finite $T_0$-space weakly equivalent to $X$ \cite{McCord1966}. Section~\ref{sec:comparison} compares the values obtained in this paper with the known sizes of minimal finite models.
\end{remark}

\begin{corollary}\label{cor:nu-product-wedge}
Let $X$ and $Y$ be spaces with $\nu(X),\nu(Y)<\infty$. Then $\nu(X\times Y)\le\nu(X)\,\nu(Y)$, and if $X$ and $Y$ are connected, $\nu(X\vee Y)\le\nu(X)+\nu(Y)-1$.
\end{corollary}

\begin{proof}
Let $\CC$ and $\DD$ be object-minimal models of $X$ and $Y$. By Lemma~\ref{lem:join}(b), $\CC\times\DD$ has $\nu(X)\nu(Y)$ objects and $B(\CC\times\DD)\cong B\CC\times B\DD\simeq X\times Y$. For the wedge, choose sources of $\CC$ and $\DD$ and apply Lemma~\ref{lem:join}(c); the homotopy type of a wedge of connected CW complexes does not depend on the vertices used as base points.
\end{proof}

The inequality for wedges is far from sharp: $\nu(S^n\vee S^n)=n+1$ by Theorem~\ref{thm:bouquets}, whereas the corollary only gives $2n+1$.

\subsection{Lower bounds from homology}\label{subsec:lower-bounds}

A nonzero class in $H_k$ needs a $k$-simplex, and hence $k+1$ objects; torsion in $H_k$ needs, in addition, a $(k+1)$-simplex whose boundary kills a multiple of a class, and hence one more object. For posets, the same dimension argument is the first step in the proof of \cite[Theorem~2.13]{BarmakMinian2007}.

\begin{proposition}\label{prop:bounds}
Let $X$ be a space and $k\ge0$.
\begin{enumerate}[(a)]
\item If $H_k(X;\Z)\neq0$, then $\nu(X)\ge k+1$.
\item If $H_k(X;\Z)$ has torsion, then $\nu(X)\ge k+2$.
\item If $X$ is connected and $n>k$, then (a) and (b) hold with $\nu_n(X)$ in place of $\nu(X)$.
\end{enumerate}
\end{proposition}

\begin{proof}
Let $\CC$ be a finite acyclic category with $H_k(B\CC;\Z)\cong H_k(X;\Z)$, and let $d=\dim B\CC$, so that $d\le|\Ob\CC|-1$ by Lemma~\ref{lem:faces}(a). This isomorphism holds if $B\CC\simeq X$, and also if there is an $n$-connected map $B\CC\to X$ with $k<n$. If $H_k(X;\Z)\neq0$, then $k\le d$. If $H_k(X;\Z)$ has torsion, then $k<d$ by Lemma~\ref{lem:faces}(c), so $|\Ob\CC|\ge d+1\ge k+2$.
\end{proof}

In the first bound for $\nu_n$ one may also take $k=n$: an $n$-connected map $B\CC\to X$ is surjective on $H_n$, so $H_n(X;\Z)\neq0$ forces $H_n(B\CC;\Z)\neq0$ and $|\Ob\CC|\ge n+1$. This is how the bound will be used for spheres and for finite groups.

\begin{remark}\label{rem:BG}
Acyclicity is essential here. For $m\ge2$ one has $\nu(B\Z/m)=\infty$, since $H_{2k-1}(B\Z/m;\Z)\cong\Z/m$ for every $k\ge1$ \cite{Hatcher2002} while $B\CC$ is finite-dimensional; thus $\nu(\RP^\infty)=\infty$, while $\nu(\RP^n)=n+1$ by Theorem~\ref{thm:RPn}. Yet the group $\Z/2=\{1,u\}$, seen as a category with one object, has classifying space $\RP^\infty$, and Lemma~\ref{lem:faces} fails for it in every part: the string $(u,\dots,u)$ is nondegenerate in every dimension, and in the normalized chain complex $\partial(u,u)=u-0+u=2u$, so torsion appears in $H_1$ with a single object. The same happens for finite monoids that are not groups; Fiedorowicz exhibited one with five elements whose classifying space is homotopy equivalent to $S^2$ \cite{Fiedorowicz2002}. For spaces such as $B\Z/m$, which have no finite model, the truncated invariants are the meaningful ones, and Section~\ref{sec:monoids} computes them for every finite group.
\end{remark}

\subsection{Stabilization}\label{subsec:stabilization}

For spaces with a finite model, the truncated invariants eventually recover $\nu$.

\begin{proposition}\label{lem:nu-n}\label{prop:stabilization}
Let $X$ be a connected space homotopy equivalent to a finite CW complex of dimension $d$. Then $\nu_n(X)=\nu(X)$ for every $n\ge\max\{d+1,\nu(X)\}$.
\end{proposition}

\begin{proof}
We may assume that $X$ is a finite CW complex of dimension $d$, and we know that $\nu_n(X)\le\nu(X)$. If $X$ is contractible, both invariants are equal to $1$. Otherwise $\nu(X)\ge2$, since the only acyclic category with one object is the terminal category, and hence $n\ge2$. Let $\CC$ be a finite acyclic category with $|\Ob\CC|=\nu_n(X)$ and an $n$-connected map $f\co B\CC\to X$. Since $n\ge2$, $f$ is an isomorphism on $\pi_1$, so it lifts to a map $\tilde f\co\widetilde{B\CC}\to\widetilde X$ between universal covers, and $\tilde f$ is again $n$-connected, since a space and its universal cover have the same homotopy groups in degrees $\ge2$.

Both universal covers have dimension at most $n-1$. Indeed, a covering space of a CW complex is a CW complex whose cells are the lifts of the cells of the base \cite{Hatcher2002}, so it has the same dimension as the base. Now $\dim\widetilde X=d\le n-1$ by hypothesis, and $\dim\widetilde{B\CC}=\dim B\CC\le|\Ob\CC|-1\le\nu(X)-1\le n-1$ by Lemma~\ref{lem:faces}(a). Hence both universal covers have zero homology in degrees $\ge n$, and since $\tilde f$ induces isomorphisms on $H_k$ for $k<n$, it is a homology isomorphism between simply connected CW complexes, hence a homotopy equivalence by Whitehead's theorem \cite[Section~4.2]{Hatcher2002}. It follows that $f$ induces isomorphisms on all homotopy groups, and by Whitehead's theorem again it is a homotopy equivalence. Therefore $\nu(X)\le|\Ob\CC|=\nu_n(X)$.
\end{proof}

The threshold is not sharp: for spheres, $\nu_n(S^k)$ already reaches $\nu(S^k)$ at $n=k$ (Remark~\ref{ex:nu-n-spheres}). For spaces without a finite model the sequence need not stabilize at all: $\nu_n(BG)=n+1$ for every nontrivial finite group $G$ (Proposition~\ref{prop:groups-minimal}).

\subsection{Object-minimality and strong minimality}\label{subsec:objmin}

Removing a beat object does not change the classifying space, so an object-minimal model cannot have one.

\begin{proposition}\label{prop:objmin}
An object-minimal model of a space is strongly minimal.
\end{proposition}

\begin{proof}
If $x$ were a beat object of an object-minimal model $\CC$ of $X$, then $B(\CC\setminus x)\simeq B\CC\simeq X$ by Lemma~\ref{lem:tanaka}(a), and $\CC\setminus x$ has fewer objects.
\end{proof}

By Lemma~\ref{lem:beat-representable}, in an object-minimal model no presheaf $\Hom(-,x)$ restricted to the objects below $x$ is representable, and no copresheaf $\Hom(x,-)$ restricted to the objects above $x$ is corepresentable; in the language of Subsection~\ref{subsec:add-sink}, no step of its construction adds a sink along a representable presheaf.

\begin{example}\label{ex:converse}
The converse of Proposition~\ref{prop:objmin} fails. The source wedge $\K\vee\K$ of two Kronecker categories, with source $0$ and sinks $1$ and $1'$, has three objects and classifying space $S^1\vee S^1$ by Lemma~\ref{lem:join}(c). Each sink receives two morphisms from $0$, neither of which factors through the other, and the morphisms out of $0$ go to two objects with no morphisms between them. By Definition~\ref{def:beat}, $\K\vee\K$ has no beat objects, so it is strongly minimal. It is not object-minimal, since $\K_3$ models $S^1\vee S^1$ with two objects (Lemma~\ref{lem:kronecker}).
\end{example}

\section{Spheres and bouquets}\label{sec:bouquets}

For posets, the minimal finite model of $S^n$ is obtained from the discrete poset with two points, which models $S^0$, by applying the non-Hausdorff suspension $n$ times \cite{BarmakMinian2007}. For acyclic categories the construction depends on the parity of the dimension: odd spheres come from the Kronecker category and its joins, while even spheres need a three-object model of $S^2$, the sphere category, which plays for $S^2$ the role that the Kronecker category plays for $S^1$. In both cases the sink lemma makes the computations transparent: the circle is the suspension of two points, and the sphere category is a cone glued to a circle by folding.

\subsection{Kronecker categories}\label{subsec:kronecker}

\begin{figure}[ht]
\centering
\begin{tikzcd}
0 \arrow[r, bend left, "f_0"] \arrow[r, bend right, "f_1"'] & 1
\end{tikzcd}
\qquad\qquad
\begin{tikzcd}[column sep=large]
0 \arrow[r, bend left=40, "f_0"] \arrow[r, phantom, "\vdots"] \arrow[r, bend right=40, "f_{r-1}"'] & 1
\end{tikzcd}

\medskip
\begin{tikzcd}[column sep=large]
0 \arrow[r, bend left, "f_0"] \arrow[r, bend right, "f_1"'] \arrow[rr, bend left=65, "h"] & 1 \arrow[r, bend left, "g_0"] \arrow[r, bend right, "g_1"'] & 2
\end{tikzcd}
\qquad\qquad
\begin{tikzcd}[column sep=large]
0 \arrow[r, bend left, "f_0"] \arrow[r, bend right, "f_1"'] \arrow[rr, bend left=65, "h"] & 1 \arrow[r, bend left=40, "g_0"] \arrow[r, phantom, "\vdots"] \arrow[r, bend right=40, "g_{r-1}"'] & 2
\end{tikzcd}
\caption{Top: the Kronecker category $\K$ and the category $\K_r$. Bottom: the sphere category $\Sph$ and the category $\Sph_r$, in which every composite $g_j\circ f_i$ is equal to $h$.}
\label{fig:kronecker}\label{fig:theta}
\end{figure}
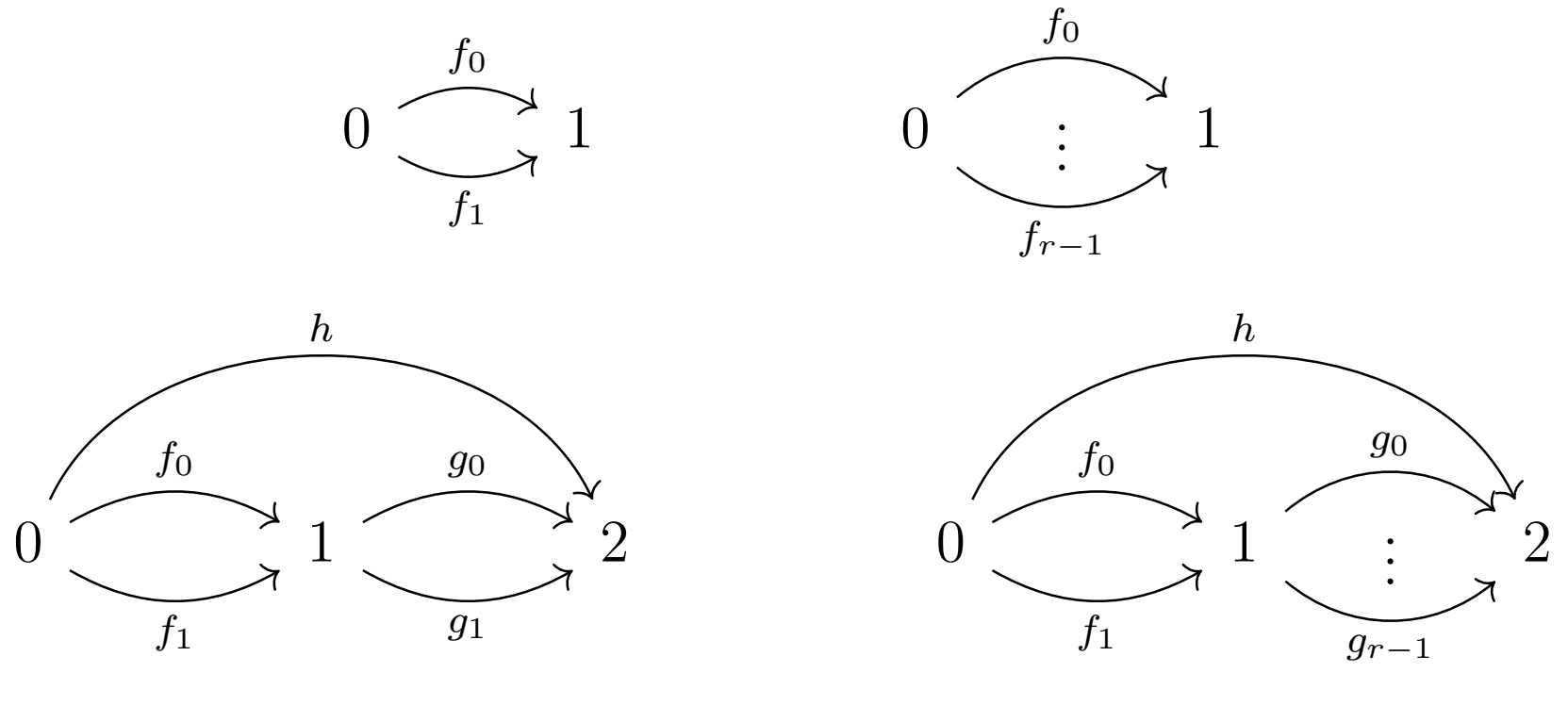

For $r\ge1$ let $\K_r$ be the category with objects $0,1$ and $r$ parallel morphisms $f_0,\dots,f_{r-1}\co0\to1$. We write $\K=\K_2$ for the \emph{Kronecker category}; see Figure~\ref{fig:kronecker}.

\begin{lemma}\label{lem:kronecker}
$B\K_r$ is a connected graph with two vertices and $r$ edges. Hence $B\K$ is homeomorphic to $S^1$, $B\K_r\simeq\bigvee^{r-1}S^1$ for $r\ge2$, and $B\K_3$ is the theta graph. Moreover $B(\K^{\star k})$ is homeomorphic to $S^{2k-1}$ for $k\ge1$.
\end{lemma}

\begin{proof}
The category $\K_r$ is obtained from the category with one object $0$ by adding the sink $1$ along the presheaf with $r$ elements, so $\DD/1$ is the discrete category on $r$ objects and $\mathrm{pr}$ collapses it to a point. By Lemma~\ref{lem:sink}, $B\K_r$ is the cone on $r$ points with its base collapsed, that is, the suspension of $r$ points, a graph with two vertices and $r$ edges. Collapsing one edge, which is a contractible subcomplex, is a homotopy equivalence onto a bouquet of $r-1$ circles. The last assertion follows from Lemma~\ref{lem:join} and $S^p*S^q\cong S^{p+q+1}$.
\end{proof}

Thus joins of Kronecker categories model odd spheres, even up to homeomorphism, with the smallest number of objects allowed by Proposition~\ref{prop:bounds}. For even spheres we need a new building block.

\subsection{The sphere category}\label{subsec:theta}

The sphere category is the analogue of $\K$ one dimension higher: adding parallel morphisms to it gives bouquets of $2$-spheres, as adding parallel morphisms to $\K$ gives bouquets of circles.

\begin{definition}\label{def:theta}
For $r\ge1$ let $\Sph_r$ be the category with objects $0,1,2$, non-trivial morphisms $f_0,f_1\co0\to1$, $g_0,\dots,g_{r-1}\co1\to2$ and $h\co0\to2$, and composites $g_j\circ f_i=h$ for all $i,j$. We write $\Sph=\Sph_2$ and call it the \emph{sphere category}. See Figure~\ref{fig:theta}.
\end{definition}

\begin{proposition}\label{prop:theta}
For $r\ge1$, $B\Sph_r\simeq\bigvee^{r-1}S^2$. In particular $B\Sph\simeq S^2$.
\end{proposition}

\begin{proof}
In the notation of Subsection~\ref{subsec:three}, $\Sph_r$ has multiplicities $(2,r)$, $c=1$, and a table with all entries $0$. Its graph $\Gamma_T$ has the vertex $h$, the vertices $g_0,\dots,g_{r-1}$, and two edges from $h$ to each $g_j$, lying over $f_0$ and $f_1$; see Figure~\ref{fig:graphs}(a) for $r=2$. It is the wedge at $h$ of $r$ bigons, each of which is mapped with degree one onto the circle $B\K$. Hence $n=r$ and $d=1$ in Proposition~\ref{prop:three}, and $B\Sph_r\simeq\bigvee^{r-1}S^2$.
\end{proof}

The sink lemma also shows what $B\Sph$ looks like: it is the cone on a figure eight, glued to the circle $B\K$ by folding both loops onto it. Four triangles meet along the edge $h$, so $B\Sph$ is not homeomorphic to $S^2$, whereas $B(\mathbf 2\star\K)\cong S^0*S^1\cong S^2$ by Lemma~\ref{lem:join}(a). Homeomorphism costs one point.

\subsection{Minimal models of spheres and bouquets}\label{subsec:bouquets}

Joining the two building blocks with Kronecker categories now gives all spheres and all bouquets of spheres of the same dimension, with the smallest possible number of objects.

\begin{theorem}\label{thm:bouquets}
Let $n,r\ge1$.
\begin{enumerate}[(a)]
\item If $n=2k+1$, then $\K_{r+1}\star\K^{\star k}$ has $n+1$ objects and classifying space homotopy equivalent to $\bigvee^rS^n$.
\item If $n=2k+2$, then $\Sph_{r+1}\star\K^{\star k}$ has $n+1$ objects and classifying space homotopy equivalent to $\bigvee^rS^n$.
\end{enumerate}
Consequently $\nu\bigl(\bigvee^rS^n\bigr)=n+1$, and in particular $\nu(S^n)=n+1$.
\end{theorem}

\begin{proof}
The object counts are $2+2k=n+1$ and $3+2k=n+1$. For $k=0$ the statements are Lemma~\ref{lem:kronecker} and Proposition~\ref{prop:theta}. For $k\ge1$, Lemmas~\ref{lem:join} and~\ref{lem:kronecker} and the facts listed after Lemma~\ref{lem:join} give
\begin{align*}
B(\K_{r+1}\star\K^{\star k})&\cong B\K_{r+1}*S^{2k-1}\simeq\Bigl(\bigvee^rS^1\Bigr)*S^{2k-1}\\
&\cong\Sigma^{2k}\Bigl(\bigvee^rS^1\Bigr)\simeq\bigvee^rS^{2k+1},
\end{align*}
and in the same way $B(\Sph_{r+1}\star\K^{\star k})\simeq\Sigma^{2k}\bigl(\bigvee^rS^2\bigr)\simeq\bigvee^rS^{2k+2}$, using Proposition~\ref{prop:theta}. The lower bound $\nu\ge n+1$ is Proposition~\ref{prop:bounds}(a).
\end{proof}

By Proposition~\ref{prop:objmin}, these categories are strongly minimal. All the topology of the bouquet is carried by a single factor of the join, $\K_{r+1}$ or $\Sph_{r+1}$: the number of objects does not grow with $r$, only the number of parallel morphisms does, whereas iterated source wedges would need $rn+1$ objects. This is the first instance of a trade-off that runs through the paper, in which objects are saved at the price of morphisms.

\begin{remark}\label{ex:nu-n-spheres}
The truncated invariants of spheres jump at once. Let $k\ge1$. If $n<k$, the inclusion of a point in $S^k$ is $n$-connected, so $\nu_n(S^k)=1$. If $n\ge k$, an $n$-connected map $B\CC\to S^k$ is surjective on $H_k$, so $\CC$ has at least $k+1$ objects, as observed after Proposition~\ref{prop:bounds}. By Theorem~\ref{thm:bouquets}, $\nu_n(S^k)=k+1=\nu(S^k)$ for $n\ge k$, one step earlier than Proposition~\ref{prop:stabilization} guarantees.
\end{remark}

\subsection{Infinitely many strongly minimal models of \texorpdfstring{$S^2$}{S2}}\label{sec:S2}

The models are not unique, and non-uniqueness already appears in the smallest case. Corollary~\ref{cor:three-S2} turns the search for three-object models of $S^2$ into a combinatorial problem about composition tables, and this problem has infinitely many solutions.

\begin{definition}\label{def:Lambda}
For $b\ge2$ let $\tSph_b$ be the three-object category with multiplicities $(2,b)$, $c=b-1$, and table with first row $0,\dots,0$ and second row $0,0,1,\dots,b-2$, that is, $T(0,j)=0$, $T(1,0)=0$ and $T(1,j)=j-1$ for $1\le j\le b-1$.
\end{definition}

For $b=2$ both rows are $0,0$, and $\tSph_2=\Sph$. In general $\tSph_b$ has the same morphisms $0\to1$ and $1\to2$ as $\Sph_b$, but the morphism $h$ of $\Sph_b$ is split into $b-1$ morphisms $h_0,\dots,h_{b-2}$; sending every $h_k$ to $h$ defines a functor $\tSph_b\to\Sph_b$, which explains the notation. Splitting $h$ changes a bouquet of $b-1$ spheres into a single sphere.

\begin{proposition}\label{prop:Lambda}
For every $b\ge2$, the category $\tSph_b$ is an object-minimal and strongly minimal model of $S^2$, and $\tSph_b$ and $\tSph_{b'}$ are not isomorphic for $b\neq b'$.
\end{proposition}

\begin{proof}
In the graph $\Gamma_T$ of $\tSph_b$, every vertex $g_j$ is joined to $h_0$ by the edge $e_{0j}$, and every vertex $h_{j-1}$ with $j\ge1$ is joined to $g_j$ by $e_{1j}$, so $\Gamma_T$ is connected; see Figure~\ref{fig:graphs}(b) for $b=3$. The edges $e_{00}$ and $e_{10}$ both join $h_0$ to $g_0$, and they lie over $f_0$ and $f_1$, so the bigon they form is mapped with degree one onto $B\K$, and $(p_T)_*$ is surjective. Since $c=b-1$, Corollary~\ref{cor:three-S2} gives $B\tSph_b\simeq S^2$. By Theorem~\ref{thm:bouquets}, $\nu(S^2)=3$, so $\tSph_b$ is object-minimal, and it is strongly minimal by Proposition~\ref{prop:objmin}. Finally, an isomorphism preserves the order $\preceq$, which is a chain on three objects, and hence the number $b$ of morphisms $1\to2$.
\end{proof}

\begin{corollary}\label{cor:infinite}
There are infinitely many isomorphism classes of object-minimal models of $S^2$, all of them strongly minimal and with three objects. They represent infinitely many strong homotopy types.
\end{corollary}

\begin{proof}
The categories $\tSph_b$ are pairwise non-isomorphic and strongly minimal, so by Lemma~\ref{lem:tanaka}(b) they are pairwise not strongly homotopy equivalent.
\end{proof}

By Lemma~\ref{lem:tanaka}(b), the strong homotopy type of a strongly minimal category determines it up to isomorphism, whereas Corollary~\ref{cor:infinite} shows that the homotopy type of its classifying space does not. Joining with Kronecker categories carries this to higher dimensions: the categories $\tSph_b\star\K^{\star k}$ are pairwise non-isomorphic object-minimal models of $S^{2k+2}$, since an isomorphism preserves the number $b$ of morphisms between the second and the third object.

\section{Moore spaces}\label{sec:moore}

Spheres are the Moore spaces of $\Z$, and the next step is to allow torsion. Moore spaces have a single non-trivial reduced homology group, but they are far from Eilenberg--MacLane spaces: $K(\Z/m,1)=B\Z/m$ has no finite model at all (Remark~\ref{rem:BG}). By Proposition~\ref{prop:bounds}, torsion costs at least one more object than a free class in the same degree, and we show that one more object is enough. In degrees one and two the models are obtained by adding a sink to a three-object category, and the sink lemma shows that the attaching map multiplies the relevant homology group by $m$; joins with Kronecker categories then give all the other degrees.

We write $\M(\Z/m,n)=S^n\cup_me^{n+1}$, where the $(n+1)$-cell is attached along a map of degree $m$. For $n\ge2$ this is the Moore space, characterized up to homotopy as the simply connected CW complex with reduced homology $\Z/m$ in degree $n$ and zero elsewhere, by the Hurewicz and Whitehead theorems \cite{Hatcher2002}. For $n=1$ we use this cell structure as the definition. By the facts listed after Lemma~\ref{lem:join}, and since the $j$-fold suspension of a map $S^n\to S^n$ of degree $m$ has degree $m$, we have $\Sigma^j\M(\Z/m,n)\simeq\M(\Z/m,n+j)$ for all $n,j\ge1$.

\subsection{Degree one}\label{subsec:moore-degree-one}

\begin{definition}\label{def:Cm1}
For $m\ge2$ let $\CC_{m,1}$ be the category with objects $0,1,2$ and non-trivial morphisms
\[
f_0,f_1\co0\to1,\qquad g_k\co1\to2,\qquad h_k\co0\to2\qquad(k\in\Z/m),
\]
with composites $g_k\circ f_i=h_{k+i}$ for $i\in\{0,1\}$ and $k\in\Z/m$.
\end{definition}

\begin{figure}[ht]
\centering
\begin{minipage}[b]{0.55\textwidth}
\centering
\begin{tikzcd}[column sep=large]
0 \arrow[r, bend left, "f_0"] \arrow[r, bend right, "f_1"'] \arrow[rr, bend left=60, "h_0" near start] \arrow[rr, bend left=85, "h_1" near end] \arrow[rr, bend right=65, "h_2"'] & 1 \arrow[r, bend left=40, "g_0"] \arrow[r, "g_1" description] \arrow[r, bend right=40, "g_2"'] & 2
\end{tikzcd}

$T=\begin{pmatrix}0&1&2\\1&2&0\end{pmatrix}$

\smallskip
(a) $\CC_{3,1}$
\end{minipage}
\hfill
\begin{minipage}[b]{0.4\textwidth}
\centering
\begin{tikzpicture}[scale=1.15]
\foreach \a/\b/\st/\lab in {90/30/solid/e_{00}, 30/-30/dashed/e_{10}, -30/-90/solid/e_{01}, -90/-150/dashed/e_{11}, 210/150/solid/e_{02}, 150/90/dashed/e_{12}}{
  \draw[\st, thick] (\a:1.2) -- (\b:1.2);
  \node at ({(\a+\b)/2}:1.42) {\scriptsize $\lab$};
}
\foreach \t in {90,30,-30,-90,-150,150} \fill (\t:1.2) circle (1.4pt);
\node[above] at (90:1.2) {\small $h_0$};
\node[right] at (30:1.2) {\small $g_0$};
\node[right] at (-30:1.2) {\small $h_1$};
\node[below] at (-90:1.2) {\small $g_1$};
\node[left] at (-150:1.2) {\small $h_2$};
\node[left] at (150:1.2) {\small $g_2$};
\end{tikzpicture}

\smallskip
(b) $\Gamma_T$ for $m=3$
\end{minipage}
\caption{(a) The category $\CC_{3,1}$ and its composition table, $g_k\circ f_i=h_{T(i,k)}$. (b) Its graph $\Gamma_T$, a hexagon; solid edges lie over $f_0$ and dashed edges over $f_1$. Going once around the hexagon goes three times around the circle $B\K$.}
\label{fig:Cm1}
\end{figure}
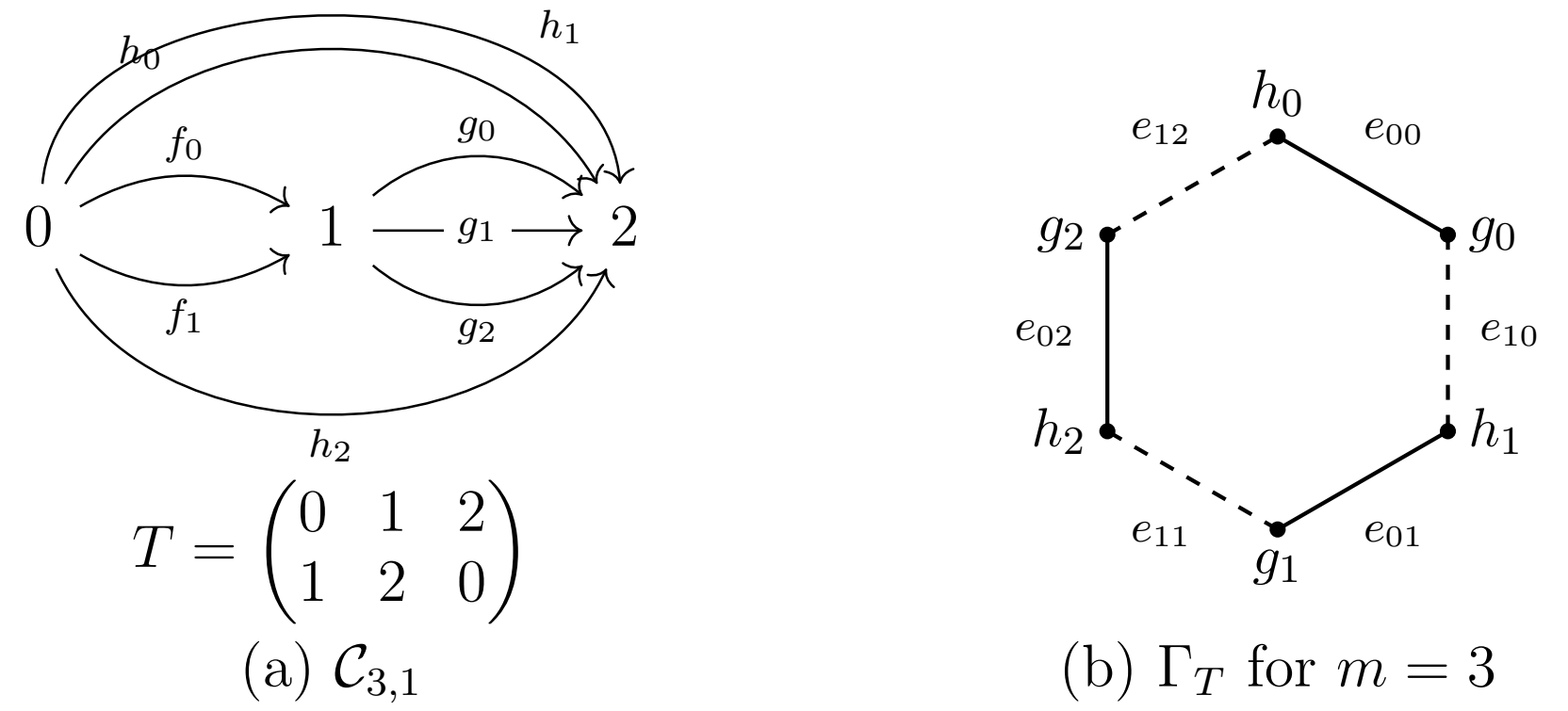

In the notation of Subsection~\ref{subsec:three}, $\CC_{m,1}$ has multiplicities $(2,m)$, $c=m$, and table $T(i,k)=i+k$ modulo $m$; see Figure~\ref{fig:Cm1}(a). It has $2m+2$ non-trivial morphisms.

\begin{theorem}\label{thm:Cm1}
For $m\ge2$, $B\CC_{m,1}$ is homeomorphic to the mapping cone of a map $S^1\to S^1$ of degree $\pm m$. Hence $B\CC_{m,1}\simeq\M(\Z/m,1)=S^1\cup_me^2$.
\end{theorem}

\begin{proof}
The graph $\Gamma_T$ has the edges $e_{0k}$ from $h_k$ to $g_k$ and $e_{1k}$ from $h_{k+1}$ to $g_k$, so it is the cycle $h_0,g_0,h_1,g_1,\dots,h_{m-1},g_{m-1}$ of length $2m$; see Figure~\ref{fig:Cm1}(b). Its fundamental class $\sum_k(e_{0k}-e_{1k})$ is mapped to $m(f_0-f_1)$, and $f_0-f_1$ generates $H_1(B\K;\Z)$, so $p_T$ is a map between circles of degree $\pm m$, and $B\CC_{m,1}$ is homeomorphic to its mapping cone by \eqref{eq:three}. Since maps $S^1\to S^1$ are classified up to homotopy by their degree, and reflecting the $2$-cell changes its sign, $B\CC_{m,1}\simeq S^1\cup_me^2$ \cite[Proposition~0.18]{Hatcher2002}.
\end{proof}

For $m=2$, the category $\CC_{2,1}$ is isomorphic to the category $\CC_{\Z/2,2}$ of Section~\ref{sec:projective}, and $B\CC_{2,1}$ is homeomorphic to $\RP^2$ by Theorem~\ref{thm:RPn}; this is the three-object model of $\RP^2$ observed by Tanaka \cite[Example~3.8]{Tanaka2018LS}.

\subsection{Degree two}\label{subsec:moore-degree-two}

In degree two the same idea needs four objects: we add a sink to a three-object category whose classifying space is $S^1\vee S^2$, and the new cone kills the circle and wraps $m$ times around the sphere.

\begin{figure}[ht]
\centering
\begin{tikzcd}[column sep=huge]
0 \arrow[r, bend left=20, "f_0"] \arrow[r, bend right=20, "f_1"'] \arrow[rr, bend left=50, "h_0"] \arrow[rr, bend right=50, "h_1"'] \arrow[rrr, bend left=75, "z"] & 1 \arrow[r, bend left=20, "g_0"] \arrow[r, bend right=20, "g_1"'] \arrow[rr, bend left=50, "y_0"] \arrow[rr, bend right=50, "y_1"'] & 2 \arrow[r, bend left=20, "x_0"] \arrow[r, bend right=20, "x_1"'] & 3
\end{tikzcd}
\caption{The category $\Q_2$. The composites are $g_j\circ f_i=h_i$, $x_k\circ g_j=y_{k+j}$, and $x_k\circ h_i=y_k\circ f_i=z$. For general $m$ there are $m$ morphisms $x_k$ and $m$ morphisms $y_k$, indexed by $k\in\Z/m$.}
\label{fig:Qm}
\end{figure}
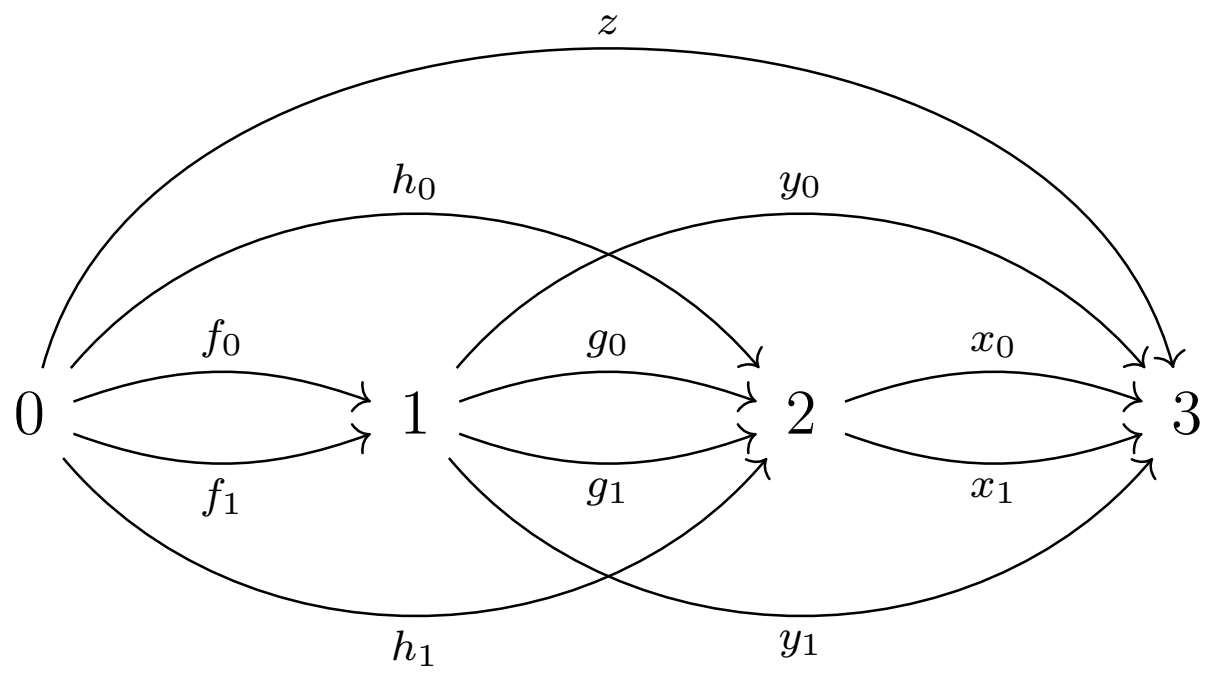

\begin{definition}\label{def:Qm}
For $m\ge2$ let $\Q_m$ be the category with objects $0,1,2,3$ and non-trivial morphisms
\begin{gather*}
f_i\co0\to1,\qquad g_j\co1\to2,\qquad h_i\co0\to2,\\
x_k\co2\to3,\qquad y_k\co1\to3,\qquad z\co0\to3,
\end{gather*}
for $i,j\in\{0,1\}$ and $k\in\Z/m$, with composites
\[
g_j\circ f_i=h_i,\qquad x_k\circ g_j=y_{k+j},\qquad x_k\circ h_i=z,\qquad y_k\circ f_i=z.
\]
See Figure~\ref{fig:Qm}.
\end{definition}

\begin{theorem}\label{thm:Qm}
For $m\ge2$, $\Q_m$ is a finite acyclic category with four objects and $2m+7$ non-trivial morphisms, and $B\Q_m\simeq\M(\Z/m,2)$.
\end{theorem}

\begin{proof}
\emph{Definition.} The objects are linearly ordered by $0<1<2<3$, the composites listed have the correct sources and targets, and associativity only has to be checked on strings $0\to1\to2\to3$, where $x_k\circ(g_j\circ f_i)=x_k\circ h_i=z=y_{k+j}\circ f_i=(x_k\circ g_j)\circ f_i$. There are $2+2+2+m+m+1=2m+7$ non-trivial morphisms. We apply Lemma~\ref{lem:sink} to the sink $t=3$.

\emph{The base category.} The full subcategory $\DD$ on $\{0,1,2\}$ has multiplicities $(2,2)$, $c=2$ and table $T(i,j)=i$. Its graph $\Gamma_T$ is a square over which $p_T$ induces zero, so $B\DD\simeq S^1\vee S^2$ by Proposition~\ref{prop:three}. We need explicit generators. With the maximal tree $\{f_0,g_0\}$, the relations $[f_i][g_j]=[h_i]$ give $[h_0]=[g_1]=1$ and $[h_1]=[f_1]$, so $\pi_1(B\DD,0)\cong\Z$ is generated by $[f_1]$ and $H_1(B\DD;\Z)$ by $f_1-f_0$. The triangles are $t_{ij}=(f_i,g_j)$, with $\partial t_{ij}=g_j-h_i+f_i$, and a $2$-chain $\sum c_{ij}t_{ij}$ is a cycle if and only if $c_{0j}+c_{1j}=0=c_{i0}+c_{i1}$; hence $H_2(B\DD;\Z)\cong\Z$ is generated by $E_0-E_1$, where $E_i=t_{i0}-t_{i1}$.

\emph{The category of elements.} The objects of $\DD/3$ are $u=(0,z)$, $p_k=(1,y_k)$ and $q_k=(2,x_k)$, and its non-trivial morphisms are $a_i^k\co u\to p_k$ over $f_i$, $c_i^k\co u\to q_k$ over $h_i$, and $b_j^k\co p_{k+j}\to q_k$ over $g_j$, for $i,j\in\{0,1\}$ and $k\in\Z/m$. The only composites are $b_j^k\circ a_i^{k+j}=c_i^k$, so the $2$-simplices are the triangles $t_{ij}^k=(a_i^{k+j},b_j^k)$, mapped by $\mathrm{pr}$ to $t_{ij}$. With the maximal tree $\{a_0^k,b_0^k\}_{k}$, the relations give $[c_0^k]=[b_1^k]=1$ and $[a_1^k]=[c_1^k]=[a_1^{k+1}]$, so $\pi_1(B(\DD/3),u)\cong\Z$ is generated by $[a_1^0]$, which is mapped to $[f_1]$: the projection induces isomorphisms on $\pi_1$ and on $H_1$. In degree two, the coefficients of $b_j^k$, $c_i^k$ and $a_i^l$ give the equations $c_{0jk}+c_{1jk}=0$, $c_{i0k}+c_{i1k}=0$ and $c_{i0l}+c_{i1,l-1}=0$ for a cycle $\sum c_{ijk}t_{ij}^k$, which force $c_{i0k}=-c_{i1k}=s_i$ with $s_1=-s_0$. Hence $H_2(B(\DD/3);\Z)\cong\Z$ is generated by $\sum_k(E_0^k-E_1^k)$, with $E_i^k=t_{i0}^k-t_{i1}^k$, and $\mathrm{pr}_*$ maps it to $m(E_0-E_1)$. Neither category has strings of three non-trivial morphisms.

\emph{Conclusion.} By Lemma~\ref{lem:sink}(d), $\pi_1(B\Q_m)\cong\pi_1(B\DD)/\langle\!\langle[f_1]\rangle\!\rangle=1$. By Lemma~\ref{lem:sink}(c), the exact sequence
\begin{gather*}
0\to\widetilde H_3(B\Q_m)\to H_2(B(\DD/3))\xrightarrow{\ \cdot m\ }H_2(B\DD)\to\widetilde H_2(B\Q_m)\to\cdots\\
\cdots\to H_1(B(\DD/3))\xrightarrow{\ \cong\ }H_1(B\DD)\to\widetilde H_1(B\Q_m)\to0
\end{gather*}
gives $H_3(B\Q_m)=0$, $H_2(B\Q_m)\cong\Z/m$ and $H_1(B\Q_m)=0$, and $H_k(B\Q_m)=0$ for $k\ge4$ by Lemma~\ref{lem:faces}(a). By the Hurewicz theorem $\pi_2(B\Q_m)\cong\Z/m$, and a based map $w\co S^2\to B\Q_m$ representing a generator extends over $S^2\cup_me^3$ to a homology isomorphism between simply connected CW complexes, which is a homotopy equivalence by Whitehead's theorem.
\end{proof}

The categorical suspension $\mathbf 2\star\CC_{m,1}$ is a model of $\M(\Z/m,2)$ with five objects; the category $\Q_m$ uses one object less.

\subsection{All degrees}\label{subsec:moore-all}

\begin{theorem}\label{thm:moore}
Let $m\ge2$ and $n\ge1$. The categories
\[
\CC_{m,1}\star\K^{\star k}\ \ (n=2k+1),\qquad \Q_m\star\K^{\star k}\ \ (n=2k+2),
\]
have $n+2$ objects and classifying spaces homotopy equivalent to $\M(\Z/m,n)$. Consequently $\nu(\M(\Z/m,n))=n+2$, and these categories are object-minimal and strongly minimal models.
\end{theorem}

\begin{proof}
The object counts are $3+2k$ and $4+2k$. For $k=0$ the statements are Theorems~\ref{thm:Cm1} and~\ref{thm:Qm}. For $k\ge1$, Lemmas~\ref{lem:join} and~\ref{lem:kronecker} give $B(\CC_{m,1}\star\K^{\star k})\cong B\CC_{m,1}*S^{2k-1}\simeq\Sigma^{2k}\M(\Z/m,1)\simeq\M(\Z/m,2k+1)$, and similarly $B(\Q_m\star\K^{\star k})\simeq\Sigma^{2k}\M(\Z/m,2)\simeq\M(\Z/m,2k+2)$. The lower bound is Proposition~\ref{prop:bounds}(b), and strong minimality follows from Proposition~\ref{prop:objmin}.
\end{proof}

\begin{remark}\label{rem:realization}
Since all these models have the object $0$ as a source, iterated source wedges (Lemma~\ref{lem:join}(c)) of the models of Theorems~\ref{thm:bouquets} and~\ref{thm:moore} realize any finite sequence $A_1,\dots,A_N$ of finitely generated abelian groups as the reduced homology of a finite acyclic category; conversely, the homology of a finite acyclic category is finitely generated and vanishes above $|\Ob\CC|-1$. Only homology groups are realized in this way: the category is not claimed to model any space given in advance, and it is far from object-minimal, since the source wedge adds objects instead of sharing them.
\end{remark}

\section{Real projective spaces and lens spaces}\label{sec:projective}

Symmetry is another way of saving points. Real projective spaces and lens spaces are quotients of spheres by free actions of finite cyclic groups, and such actions can already be seen on posets that model the spheres. The quotient of a poset by a group acting freely on its objects is an acyclic category whose classifying space is the orbit space, but it is no longer a poset: the orbits of the morphisms between two orbits of objects become parallel morphisms. The smallest example is the crown, with minimal elements $a_0,a_1$, maximal elements $b_0,b_1$ and $a_r<b_s$ for all $r,s$. It is a model of $S^1$ on which the swap $a_0\leftrightarrow a_1$, $b_0\leftrightarrow b_1$ is the antipodal map, and its four morphisms fall into the two orbits $\{a_0\to b_0,\,a_1\to b_1\}$ and $\{a_0\to b_1,\,a_1\to b_0\}$, which become the two parallel morphisms of the quotient $\K$, a model of $\RP^1$:
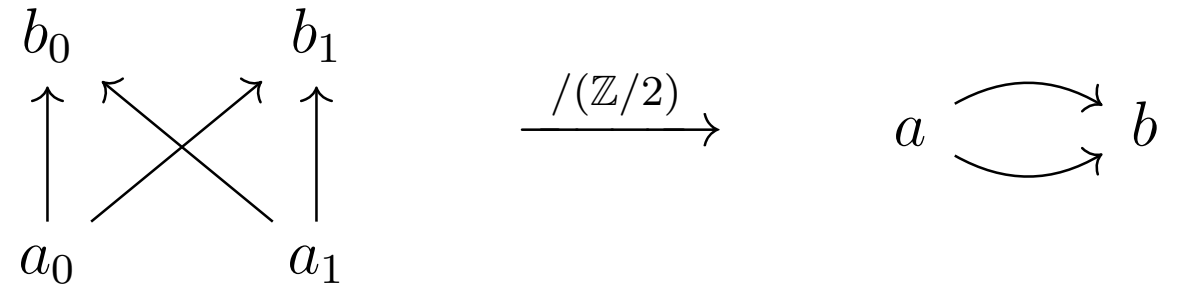
\[
\begin{tikzcd}[column sep=2.5em, row sep=2.2em]
b_0 & b_1 \\
a_0 \arrow[u] \arrow[ur] & a_1 \arrow[u] \arrow[ul]
\end{tikzcd}
\qquad\xrightarrow{\ /(\Z/2)\ }\qquad
\begin{tikzcd}
a \arrow[r, bend left] \arrow[r, bend right] & b
\end{tikzcd}
\]
In general the number of objects is divided by the order of the group, and the resulting models have one object more than the dimension, the lower bound of Proposition~\ref{prop:bounds}.

Kozlov studies quotients of acyclic categories by group actions \cite[Chapter~14]{Kozlov2008}: the quotient $\CC/G$ is the colimit of the action \cite[Definition~14.8]{Kozlov2008}, the orbit space $\nerve(\CC)/G$ is a regular trisp \cite[Proposition~14.4, (14.1)]{Kozlov2008}, and the canonical surjection $c\co\nerve(\CC)/G\to\nerve(\CC/G)$ is an isomorphism if and only if Condition~(C1) of \cite[Theorem~14.14]{Kozlov2008} holds. For actions that are free on objects, (C1) holds, since an element of $G$ that sends a morphism to another one with the same source fixes that source, and the quotient has the following explicit description \cite[Proposition~14.13]{Kozlov2008}. We sketch the argument to fix the conventions.

\begin{lemma}\label{lem:free-quotient}
Let a finite group $G$ act on a finite acyclic category $\CC$ so that the action on $\Ob\CC$ is free. Let $\CC/G$ be the category whose objects are the orbits $Gx$ of objects, whose morphisms $Gx\to Gy$ are the orbits $Gf$ of morphisms $f\co x'\to y'$ with $x'\in Gx$ and $y'\in Gy$, and whose composition is $Gg\circ Gf=G(h(g)\circ f)$, where $h\in G$ is the unique element such that the source of $h(g)$ is the target of $f$. Then $\CC/G$ is a finite acyclic category, and the map sending the orbit of a string $(f_1,\dots,f_k)$ of $\CC$ to $(Gf_1,\dots,Gf_k)$ is an isomorphism of trisps $\nerve(\CC)/G\to\nerve(\CC/G)$. Consequently $B(\CC/G)$ is homeomorphic to the orbit space $B\CC/G$.
\end{lemma}

\begin{proof}
Since $G$ acts freely on objects, for objects $y',y''$ in the same orbit there is a unique $h\in G$ with $h(y'')=y'$. Replacing $f$ and $g$ by $k(f)$ and $l(g)$ replaces $h$ by $khl^{-1}$, and $(khl^{-1})(l(g))\circ k(f)=k(h(g)\circ f)$, so the composition is well defined. A string $(Gf_1,\dots,Gf_k)$ of composable morphisms of $\CC/G$ lifts to strings $(f_1,\dots,f_k)$ of $\CC$: choose a representative $f_1$ and then, successively, the unique representative of the next orbit whose source is the target of the previous morphism. The lifts form a single orbit, permuted freely by $G$, and $Gf_{l+1}\circ Gf_l=G(f_{l+1}\circ f_l)$ for a lift, so associativity and identities follow from those of $\CC$. A non-trivial endomorphism of $\CC/G$, or a pair of non-trivial morphisms in opposite directions, would lift to a string of non-trivial morphisms from an object $x$ to $h(x)$ for some $h\in G$, and then $x\prec h(x)\prec\cdots\prec h^e(x)=x$ for the order $e$ of $h$, where $\prec$ is the strict order induced by $\preceq$; hence $\CC/G$ is acyclic. Finally, lifting identifies the $k$-simplices of $\nerve(\CC/G)$ with the $G$-orbits of $k$-simplices of $\nerve(\CC)$, compatibly with faces, and $G$ acts freely on the simplices of $\nerve(\CC)$, preserving the order of their vertices. Hence $B\CC/G$ is obtained by gluing one simplex for each orbit, and it is the realization of $\nerve(\CC)/G\cong\nerve(\CC/G)$.
\end{proof}

\subsection{Real projective spaces}\label{subsec:RPn}

For $n\ge1$ let $\CC_{\Z/2,n}$ be the category with objects $0,1,\dots,n$, with $\Hom(i,j)=\{g_{ij}^0,g_{ij}^1\}$ for $i<j$ and no non-trivial morphisms otherwise, and with composites $g_{jk}^{s'}\circ g_{ij}^s=g_{ik}^{s+s'}$, the exponents taken in $\Z/2$. Composition is associative because addition in $\Z/2$ is, and $\CC_{\Z/2,n}$ is acyclic with $n+1$ objects and $n(n+1)$ non-trivial morphisms. Note that $\CC_{\Z/2,1}\cong\K$ and $\CC_{\Z/2,2}\cong\CC_{2,1}$; in Section~\ref{sec:monoids}, $\CC_{\Z/2,n}$ reappears as the case $M=\Z/2$ of a construction that applies to every finite monoid.

\begin{figure}[ht]
\centering
\begin{minipage}[b]{0.26\textwidth}
\centering
\begin{tikzpicture}[>={Stealth[length=4pt]}, every node/.style={font=\small}]
\foreach \i in {0,1,2} {
  \node (p\i0) at (0,1.25*\i) {$(\i,0)$};
  \node (p\i1) at (1.9,1.25*\i) {$(\i,1)$};
}
\foreach \i/\j in {0/1,1/2} \foreach \a in {0,1} \foreach \b in {0,1}
  \draw[->] (p\i\a) -- (p\j\b);
\end{tikzpicture}

\smallskip
(a) $\PP_2$, a model of $S^2$
\end{minipage}
\hfill
\begin{minipage}[b]{0.33\textwidth}
\centering
\begin{tikzcd}[column sep=2.4em]
0 \arrow[r, bend left=25, "0"] \arrow[r, bend right=25, "1"'] \arrow[rr, bend left=75, "0"] \arrow[rr, bend right=75, "1"'] & 1 \arrow[r, bend left=25, "0"] \arrow[r, bend right=25, "1"'] & 2
\end{tikzcd}

\smallskip
(b) $\CC_{\Z/2,2}=\PP_2/(\Z/2)$
\end{minipage}
\hfill
\begin{minipage}[b]{0.33\textwidth}
\centering
\begin{tikzpicture}[scale=1.6, >={Stealth[length=4pt]}, lbl/.style={font=\scriptsize}]
\fill[gray!12] (1,0) -- (0,1) -- (-1,0) -- (0,-1) -- cycle;
\draw (1,0) -- (0,1) -- (-1,0) -- (0,-1) -- cycle;
\draw (1,0) -- (-1,0); \draw (0,1) -- (0,-1);
\draw[->] (1,0) -- (0.45,0.55); \node[lbl, above right] at (0.5,0.5) {$0$};
\draw[->] (-1,0) -- (-0.45,0.55); \node[lbl, above left] at (-0.5,0.5) {$1$};
\draw[->] (-1,0) -- (-0.45,-0.55); \node[lbl, below left] at (-0.5,-0.5) {$0$};
\draw[->] (1,0) -- (0.45,-0.55); \node[lbl, below right] at (0.5,-0.5) {$1$};
\draw[->] (1,0) -- (0.5,0); \node[lbl, above] at (0.7,0) {$0$};
\draw[->] (-1,0) -- (-0.5,0); \node[lbl, above] at (-0.7,0) {$1$};
\draw[->] (0,1) -- (0,0.5); \node[lbl, right] at (0,0.72) {$0$};
\draw[->] (0,-1) -- (0,-0.5); \node[lbl, right] at (0,-0.72) {$1$};
\foreach \p in {(1,0),(0,1),(-1,0),(0,-1),(0,0)} \fill \p circle (1.2pt);
\node[right] at (1,0) {\small $0$};
\node[above] at (0,1) {\small $1$};
\node[left] at (-1,0) {\small $0$};
\node[below] at (0,-1) {\small $1$};
\node[below left] at (0,0) {\small $2$};
\end{tikzpicture}

\smallskip
(c) $B\CC_{\Z/2,2}\cong\RP^2$
\end{minipage}
\caption{(a) The minimal finite model $\PP_2$ of $S^2$; the group $\Z/2$ swaps the two columns, and on the octahedron $B\PP_2$ it acts as the antipodal map. (b) The quotient $\CC_{\Z/2,2}$, Tanaka's three-object model of $\RP^2$; an arrow $i\to j$ labeled $s$ is $g_{ij}^s$. (c) Its classifying space, one hemisphere of the octahedron: the four triangles are the composable pairs $(g_{01}^s,g_{12}^{s'})$, with third edge $g_{02}^{s+s'}$, and opposite boundary edges, which carry the same label, are identified antipodally.}
\label{fig:RP2}
\end{figure}
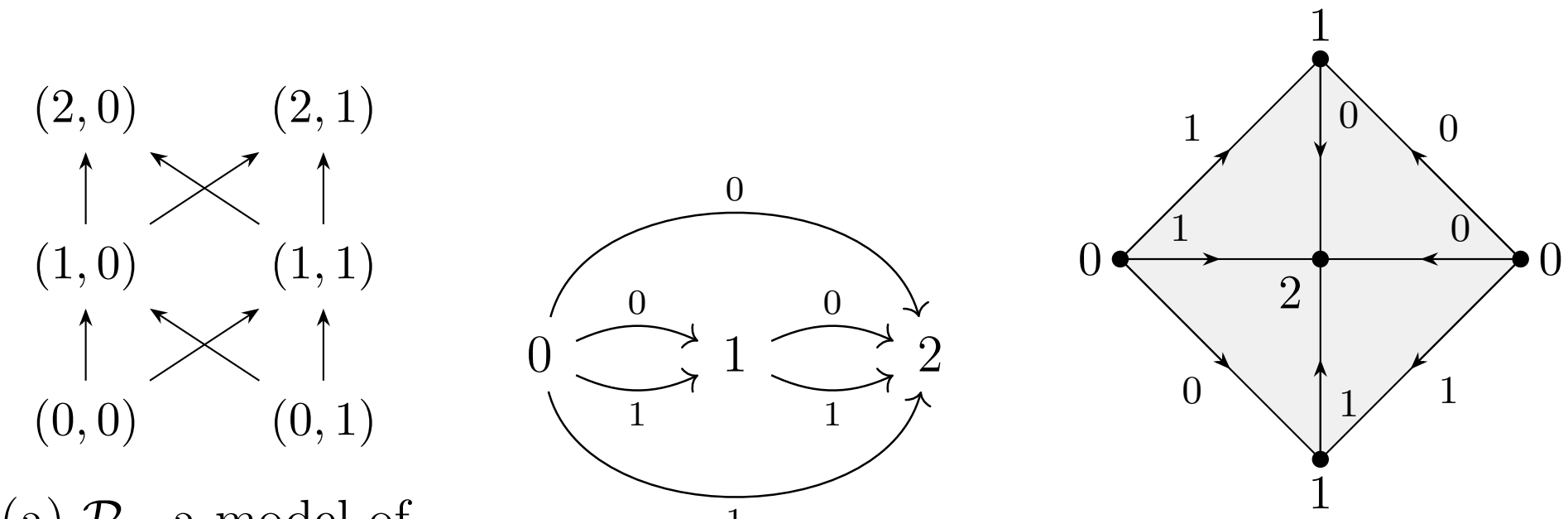

Let $\PP_n$ be the stack of $n+1$ copies of $\mathbf 2$, with objects $(i,a)$ for $0\le i\le n$ and $a\in\Z/2$, and one morphism $(i,a)\to(j,b)$ whenever $i<j$. The group $\Z/2$ acts on $\PP_n$ by $(i,a)\mapsto(i,a+1)$, simultaneously swapping the two objects of each copy of $\mathbf 2$, as in \cite[Example~14.6]{Kozlov2008}. The poset $\PP_n$ is the $n$-fold non-Hausdorff suspension of the discrete poset with two points, which is the minimal finite model of $S^n$ \cite{BarmakMinian2007}. Thus we are letting $\Z/2$ act on the minimal poset models of the spheres, and the quotient collapses their $2n+2$ points to $n+1$ objects. For $n=2$, $\PP_2$ is the six-point model of $S^2$, whose nerve is the octahedron, and the quotient $\CC_{\Z/2,2}$ is the three-object model of $\RP^2$ observed by Tanaka \cite[Example~3.8]{Tanaka2018LS}; see Figure~\ref{fig:RP2}.

\begin{theorem}\label{thm:RPn}
For every $n\ge1$, $B\CC_{\Z/2,n}$ is homeomorphic to $\RP^n$. Consequently $\nu(\RP^n)=n+1$, and $\CC_{\Z/2,n}$ is an object-minimal and strongly minimal model of $\RP^n$.
\end{theorem}

\begin{proof}
The action is free on objects. The orbit of the morphism $(i,a)\to(j,b)$ is determined by $i<j$ and by $a+b\in\Z/2$, and if $(i,a)\to(j,b)$ and $(j,b)\to(l,c)$ are composable, then $(a+b)+(b+c)=a+c$. Hence $\PP_n/(\Z/2)$, in the sense of Lemma~\ref{lem:free-quotient}, is isomorphic to $\CC_{\Z/2,n}$, and $B\CC_{\Z/2,n}$ is homeomorphic to $B\PP_n/(\Z/2)$.

Since $\PP_n$ is a poset, $\nerve(\PP_n)$ is the simplicial complex whose simplices are the sets of objects with pairwise distinct first coordinates. Sending $(i,a)$ to $(-1)^ae_i\in\mathbb R^{n+1}$ and extending linearly on simplices identifies $B\PP_n$ with the boundary of the cross-polytope with vertices $\pm e_0,\dots,\pm e_n$, and hence with $S^n$ by radial projection. Under this identification the action is $v\mapsto-v$, so it is antipodal, as observed in \cite[Example~14.6]{Kozlov2008}. Therefore $B\CC_{\Z/2,n}$ is homeomorphic to $S^n/\{\pm1\}=\RP^n$.

For the lower bound, $H_n(\RP^n;\Z)\cong\Z$ if $n$ is odd and $H_{n-1}(\RP^n;\Z)\cong\Z/2$ if $n\ge2$ is even \cite{Hatcher2002}, so $\nu(\RP^n)\ge n+1$ by Proposition~\ref{prop:bounds}. Strong minimality follows from Proposition~\ref{prop:objmin}.
\end{proof}

\begin{remark}\label{rem:RP-cells}
The sink lemma shows how these models grow. The full subcategory of $\CC_{\Z/2,n}$ on $\{0,\dots,n-1\}$ is $\CC_{\Z/2,n-1}$, and $n$ is a sink. Write $(i,s)$ for the object $(i,g_{in}^s)$ of $\DD/n$. A morphism $(i,s)\to(j,s')$ is a morphism $g_{ij}^a$ with $g_{jn}^{s'}\circ g_{ij}^a=g_{in}^s$, that is, with $a=s+s'$, so there is exactly one when $i<j$ and none otherwise. Hence $\DD/n$ is the poset $\PP_{n-1}$, and $\mathrm{pr}$ is the quotient functor $\PP_{n-1}\to\PP_{n-1}/(\Z/2)\cong\CC_{\Z/2,n-1}$. By the sink lemma, $B\CC_{\Z/2,n}$ is obtained from $\RP^{n-1}$ by attaching the cone on $S^{n-1}$ along the double covering $S^{n-1}\to\RP^{n-1}$, which is the classical way of attaching the top cell of $\RP^n$ \cite[Example~0.4]{Hatcher2002}. Each new object adds one cell.
\end{remark}

\subsection{Lens spaces}\label{subsec:lens}

Let $m\ge2$ and $k\ge0$, let $q_0,\dots,q_k$ be integers prime to $m$, and let $\zeta=\exp(2\pi\sqrt{-1}/m)$. The lens space $L(m;q_0,\dots,q_k)$ is the quotient of the unit sphere $S^{2k+1}\subset\mathbb C^{k+1}$ by the free action of $\Z/m$ whose generator acts by $(z_0,\dots,z_k)\mapsto(\zeta^{q_0}z_0,\dots,\zeta^{q_k}z_k)$. It is a closed orientable manifold of dimension $2k+1$; for $m=2$ it is $\RP^{2k+1}$, and for $k=0$ it is a circle. For $k\ge1$ its fundamental group is $\Z/m$, its universal cover is $S^{2k+1}$, and $H_i(L(m;q_0,\dots,q_k);\Z)$ is $\Z$ for $i=0,2k+1$, $\Z/m$ for odd $0<i<2k+1$, and zero otherwise \cite[Example~2.43]{Hatcher2002}. These invariants depend only on $m$ and $k$, but the homotopy type and the homeomorphism type of $L(m;q_0,\dots,q_k)$ also depend on the parameters $q_i$ \cite{Whitehead1941,Reidemeister1935,Brody1960}. This is why our models retain the whole action, and not only the group.

For $m\ge2$ let $\mathcal O_m$ be the \emph{crown} of size $m$: the poset with $2m$ elements $a_r,b_r$ for $r\in\Z/m$, with $a_r<b_r$ and $a_{r-1}<b_r$. For $m=2$ it is the four-point crown above, the minimal finite model of $S^1$. Its nerve is a cycle with $2m$ vertices, in which $b_r$ is joined to $a_{r-1}$ and to $a_r$, so $B\mathcal O_m$ is a circle. For $q$ prime to $m$ we write $\bar q\in\Z/m$ for the inverse of $q$ modulo $m$.

\begin{figure}[ht]
\centering
\begin{tikzpicture}[scale=0.9, >={Stealth[length=5pt]},
  aob/.style={circle, fill=black, inner sep=1.5pt},
  bob/.style={circle, draw=black, fill=white, inner sep=1.5pt},
  lbl/.style={font=\scriptsize}]
\foreach \i/\y in {0/0,1/3.2} {
  \begin{scope}[yshift=\y cm]
  \draw[thick] (0,0) -- (-0.5,1.1); \draw[thick, dashed] (0,0) -- (0.5,1.1);
  \draw[thick] (1,0) -- (0.5,1.1);  \draw[thick, dashed] (1,0) -- (1.5,1.1);
  \draw[thick] (3,0) -- (2.5,1.1);  \draw[thick, dashed] (3,0) -- (3.5,1.1);
  \draw[thick] (1.5,1.1) -- (1.72,0.62); \draw[thick, dashed] (2.5,1.1) -- (2.28,0.62);
  \node[aob] at (0,0) {}; \node[aob] at (1,0) {}; \node[aob] at (3,0) {};
  \node[bob] at (-0.5,1.1) {}; \node[bob] at (0.5,1.1) {}; \node[bob] at (1.5,1.1) {};
  \node[bob] at (2.5,1.1) {}; \node[bob, draw=gray] at (3.5,1.1) {};
  \node at (2,0) {$\cdots$};
  \node[lbl, below] at (0,0) {$a_0$}; \node[lbl, below] at (1,0) {$a_1$}; \node[lbl, below] at (3,0) {$a_{m-1}$};
  \node[lbl, above] at (-0.5,1.1) {$b_0$}; \node[lbl, above] at (0.5,1.1) {$b_1$}; \node[lbl, above] at (1.5,1.1) {$b_2$};
  \node[lbl, above] at (2.5,1.1) {$b_{m-1}$}; \node[lbl, above, gray] at (3.5,1.1) {$b_0$};
  \node[left] at (-0.9,0.55) {$\mathcal O_m^{(\i)}$};
  \node[lbl, gray, right] at (3.9,0.55) {$r\mapsto r+q_{\i}$};
  \end{scope}
}
\draw[->, double, thick, gray] (1.5,1.85) -- node[lbl, right=2pt, black, align=left]{all relations $x<y$,\\ $x$ below, $y$ above} (1.5,2.75);
\node at (1.5,-1.0) {(a) $\PP$, for $k=1$};
\draw[->, gray, thick] (5.55,2.15) -- node[above, lbl, black]{$/(\Z/m)$} (6.75,2.15);
\begin{scope}[xshift=8.1cm]
\node[aob, label=left:{$A_0$}] (A0) at (0,0) {};
\node[bob, label=left:{$B_0$}] (B0) at (0,1.1) {};
\node[aob, label=left:{$A_1$}] (A1) at (0,3.2) {};
\node[bob, label=left:{$B_1$}] (B1) at (0,4.3) {};
\draw[->, thick] (A0) to[bend left=35] node[lbl, left]{$0$} (B0);
\draw[->, thick, dashed] (A0) to[bend right=35] node[lbl, right]{$\bar q_0$} (B0);
\draw[->, thick] (A1) to[bend left=35] node[lbl, left]{$0$} (B1);
\draw[->, thick, dashed] (A1) to[bend right=35] node[lbl, right]{$\bar q_1$} (B1);
\draw[->, double, thick, gray] (0,1.55) -- node[lbl, right, black, align=left]{$m$ morphisms $X\to Y$,\\ labels in $\Z/m$} (0,2.75);
\node at (0.3,-1.0) {(b) $\mathcal L_m(q_0,q_1)$};
\end{scope}
\end{tikzpicture}
\caption{(a) The stack $\PP$ of the proof of Theorem~\ref{thm:lens} for $k=1$: two crowns $\mathcal O_m$, each drawn unrolled with $b_0$ repeated, in which every element of the lower copy is below every element of the upper one (these relations are not drawn). Solid edges are $a_r\to b_r$ and dashed edges are $a_r\to b_{r+1}$. The generator of $\Z/m$ acts by $r\mapsto r+q_i$ on the $i$-th copy, and $B\PP\cong S^1*S^1=S^3$. (b) The quotient $\mathcal L_m(q_0,q_1)$. Each row of (a) becomes one object, and the solid and dashed edges of the $i$-th copy become the two morphisms $A_i\to B_i$, labeled by the phase differences $0$ and $\bar q_i$. The $m^2$ relations between two rows of different copies form $m$ orbits, labeled by $\Z/m$. For $k=0$ only one crown remains, and its quotient is the Kronecker category $\K$.}
\label{fig:lens-stack}
\end{figure}
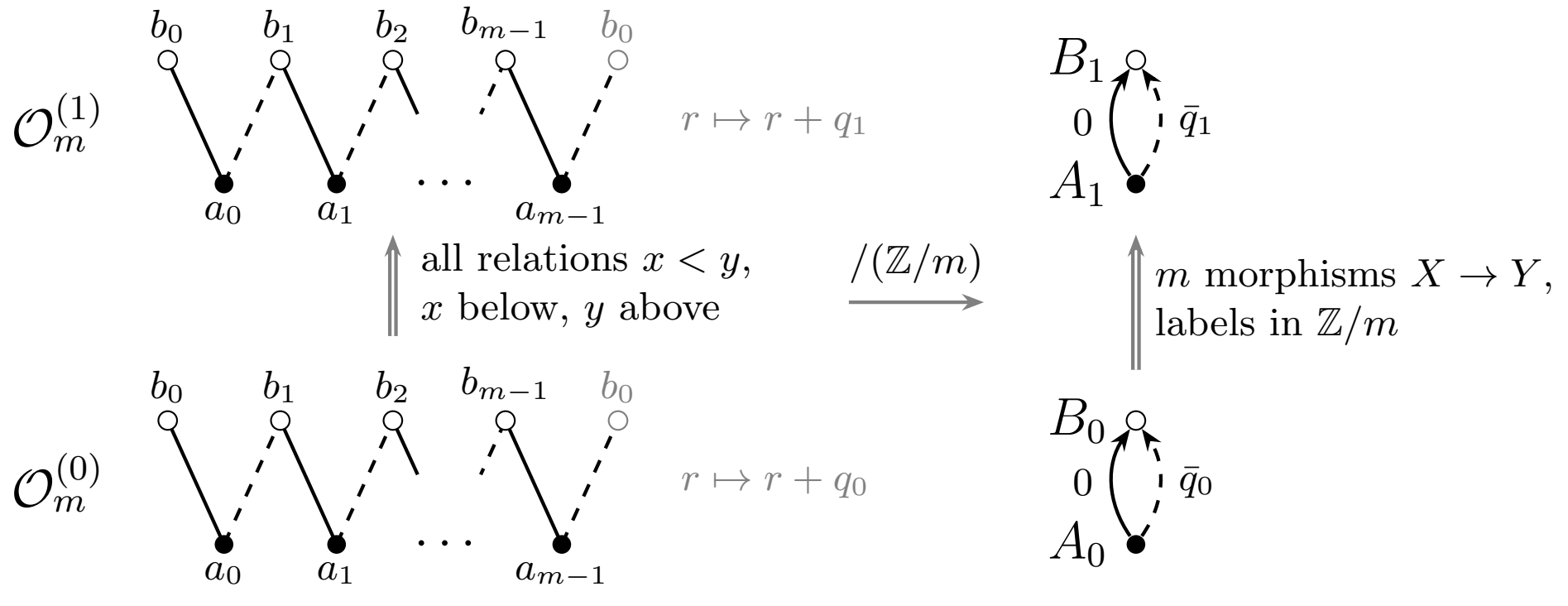

\begin{definition}\label{def:lens}
Let $\mathcal L=\mathcal L_m(q_0,\dots,q_k)$ be the category with objects $A_0,B_0,A_1,B_1,\dots,A_k,B_k$, whose non-trivial morphisms are labeled by elements of $\Z/m$ as follows: the morphisms $A_i\to B_i$ are labeled by $0$ and $\bar q_i$; for $i<j$ and $X\in\{A_i,B_i\}$, $Y\in\{A_j,B_j\}$, the morphisms $X\to Y$ are labeled by all elements of $\Z/m$; there are no other non-trivial morphisms. The composite of morphisms labeled $d$ and $d'$ is the morphism labeled $d+d'$. See Figure~\ref{fig:lens-stack}(b).
\end{definition}

The composite of two non-trivial morphisms $X\to Y\to Z$ goes between objects of different indices, since there is no object strictly between $A_i$ and $B_i$, so it is defined; associativity follows from that of addition. The objects are linearly ordered by $A_0<B_0<\cdots<A_k<B_k$, so $\mathcal L$ is acyclic. It has $2k+2$ objects and $2(k+1)+4m\binom{k+1}{2}=2(k+1)(mk+1)$ non-trivial morphisms. For $k=0$ it is isomorphic to $\K$, and for $m=2$ it is isomorphic to $\CC_{\Z/2,2k+1}$.

\begin{theorem}\label{thm:lens}
Let $m\ge2$, $k\ge0$, and let $q_0,\dots,q_k$ be integers prime to $m$. Then $B\mathcal L_m(q_0,\dots,q_k)$ is homeomorphic to $L(m;q_0,\dots,q_k)$. Consequently $\nu(L)=2k+2$ for every lens space $L$ of dimension $2k+1$, and $\mathcal L_m(q_0,\dots,q_k)$ is an object-minimal and strongly minimal model of $L(m;q_0,\dots,q_k)$.
\end{theorem}

\begin{proof}
Let $\PP$ be the stack of $k+1$ copies of $\mathcal O_m$ (Figure~\ref{fig:lens-stack}(a) for $k=1$), with objects $a_r^{(i)},b_r^{(i)}$ for $0\le i\le k$ and $r\in\Z/m$. Let the generator of $\Z/m$ act by $r\mapsto r+q_i$ on the objects of the $i$-th copy. This is an automorphism of $\PP$, because it preserves each copy and its order, and the morphisms between different copies are unique. It is free on objects, since $q_i$ is prime to $m$.

Give the objects $a_r^{(i)}$ and $b_r^{(i)}$ the \emph{phase} $\bar q_ir\in\Z/m$. The generator adds $1$ to every phase. Hence the orbit of a morphism $x\to y$ of $\PP$ is determined by the orbits of $x$ and $y$ and by the difference of phases of $y$ and $x$, and differences of phases add under composition. For morphisms between different copies all differences occur, and the morphisms $a_r^{(i)}\to b_r^{(i)}$ and $a_{r-1}^{(i)}\to b_r^{(i)}$ have differences $0$ and $\bar q_i$. Sending the orbit of $a_r^{(i)}$ to $A_i$, that of $b_r^{(i)}$ to $B_i$, and the orbit of a morphism to its difference of phases, we obtain an isomorphism $\PP/(\Z/m)\cong\mathcal L_m(q_0,\dots,q_k)$, where the quotient is that of Lemma~\ref{lem:free-quotient}. By that lemma, $B\mathcal L$ is homeomorphic to $B\PP/(\Z/m)$.

Realize $B\mathcal O_m$, in the $i$-th copy, as the regular $2m$-gon in $\mathbb C$ with $a_r$ at $\zeta^r$ and $b_r$ at $\exp\bigl(\pi\sqrt{-1}(2r-1)/m\bigr)$, between $a_{r-1}$ and $a_r$. The generator then acts as multiplication by $\zeta^{q_i}$, and radial projection is an equivariant homeomorphism onto the unit circle with this rotation. By Lemma~\ref{lem:join}, $B\PP$ is homeomorphic to the join of the $k+1$ polygons, and the homeomorphism is defined simplex by simplex, so it is equivariant for the diagonal action. Joining the radial projections gives an equivariant homeomorphism onto the join of $k+1$ unit circles. Finally, the map from the join of $k+1$ unit circles to $S^{2k+1}$ sending $\sum_it_iz_i$, with $t_i\ge0$ and $\sum_it_i=1$, to $(\sqrt{t_0}\,z_0,\dots,\sqrt{t_k}\,z_k)$ is a homeomorphism, equivariant for the diagonal rotations. Hence $B\PP$ is equivariantly homeomorphic to $S^{2k+1}$ with the action defining $L(m;q_0,\dots,q_k)$, and $B\mathcal L\cong L(m;q_0,\dots,q_k)$.

Since $H_{2k+1}(L;\Z)\cong\Z$, Proposition~\ref{prop:bounds} gives $\nu(L)\ge2k+2$ for every lens space $L$ of dimension $2k+1$, and every such $L$ is of the form $L(m;q_0,\dots,q_k)$. Strong minimality follows from Proposition~\ref{prop:objmin}.
\end{proof}

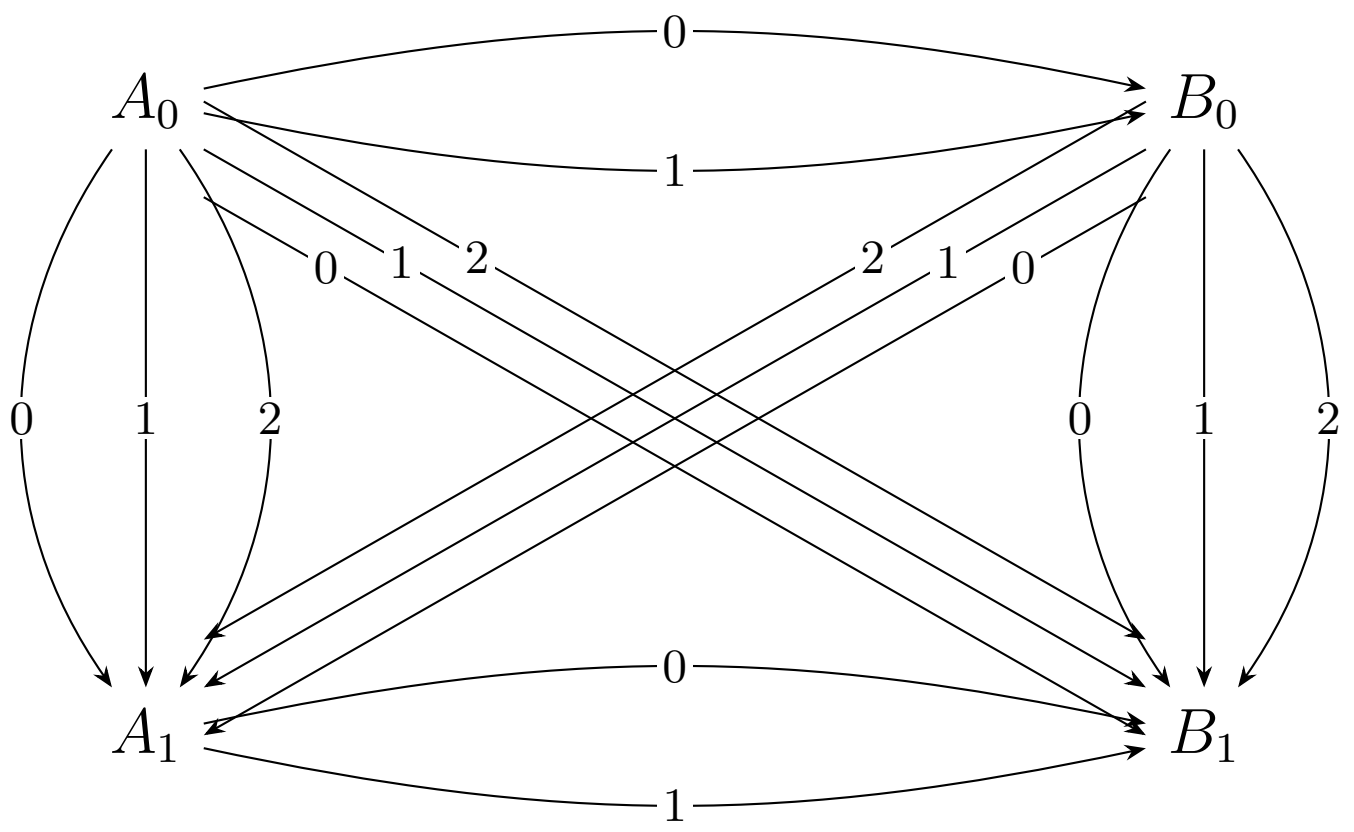
\begin{figure}[ht]
\centering
\begin{tikzpicture}[>={Stealth[length=4pt]}, lbl/.style={font=\scriptsize, fill=white, inner sep=1pt}]
\node (A0) at (0,4.2) {$A_0$};
\node (B0) at (7,4.2) {$B_0$};
\node (A1) at (0,0) {$A_1$};
\node (B1) at (7,0) {$B_1$};
\draw[->] (A0) to[bend left=12] node[lbl]{$0$} (B0);
\draw[->] (A0) to[bend right=12] node[lbl]{$1$} (B0);
\draw[->] (A1) to[bend left=12] node[lbl]{$0$} (B1);
\draw[->] (A1) to[bend right=12] node[lbl]{$1$} (B1);
\draw[->] (A0) to[bend right=35] node[lbl]{$0$} (A1);
\draw[->] (A0) -- node[lbl]{$1$} (A1);
\draw[->] (A0) to[bend left=35] node[lbl]{$2$} (A1);
\draw[->] (B0) to[bend right=35] node[lbl]{$0$} (B1);
\draw[->] (B0) -- node[lbl]{$1$} (B1);
\draw[->] (B0) to[bend left=35] node[lbl]{$2$} (B1);
\foreach \s/\l/\p in {-9pt/0/0.13,0pt/1/0.21,9pt/2/0.29} {
  \draw[->] ([yshift=\s]A0.south east) -- node[lbl, pos=\p]{$\l$} ([yshift=\s]B1.north west);
  \draw[->] ([yshift=\s]B0.south west) -- node[lbl, pos=\p]{$\l$} ([yshift=\s]A1.north east);
}
\end{tikzpicture}
\caption{The category $\mathcal L_3(1,1)$, an object-minimal model of the lens space $L(3;1,1)$, with all its $16$ non-trivial morphisms. Each morphism is labeled by an element of $\Z/3$, and labels add under composition. The full subcategory on $A_0$, $B_0$ and $A_1$ is $\CC_{3,1}$, a model of $\M(\Z/3,1)$.}
\label{fig:L31}
\end{figure}

As for projective spaces, the sink lemma recovers the minimal cell structure. In $\mathcal L_m(1,1)$, for instance, the full subcategory on $A_0,B_0,A_1$ is $\CC_{m,1}$, with $f_i$, $g_k$ and $h_k$ the morphisms $A_0\to B_0$, $B_0\to A_1$ and $A_0\to A_1$ labeled $i$, $k$ and $k$, so these three objects produce the $2$-skeleton $S^1\cup_me^2$. The object $B_1$ is a sink, and a morphism $(X,u)\to(Y,v)$ of $\DD/B_1$ is the morphism $X\to Y$ labeled $u-v$, when it exists. Hence $\DD/B_1$ is a poset: the objects $(A_0,u)$ and $(B_0,v)$ form a crown isomorphic to $\mathcal O_m$, with $(A_0,u)<(B_0,v)$ if and only if $u-v\in\{0,1\}$, and the two objects $(A_1,0)$ and $(A_1,1)$ lie above all of them. Its classifying space is the suspension of a circle, and by Lemma~\ref{lem:sink} the last object attaches the $3$-cell. Figure~\ref{fig:L31} shows the case $m=3$.

\section{Finite approximations to classifying spaces}\label{sec:monoids}

Classifying spaces of finite groups have no finite acyclic models (Remark~\ref{rem:BG}), and exact minimization says nothing about them. Truncated minimization does. In this section we spread a finite monoid $M$ over $n+1$ objects and obtain a finite acyclic category with an $n$-connected map to $BM$; for a finite group $G$ this is the $n$-th stage of Milnor's construction of $BG$ \cite{Milnor1956}, and it realizes $\nu_n(BG)$. The construction is a strict form of subdivision, and the proofs only use the sink lemma.

\subsection{Spreading a monoid over \texorpdfstring{$n+1$}{n+1} objects}\label{subsec:spreading}

\begin{definition}\label{def:CMn}
Let $X$ be a small category and let $n\in\mathbb N\cup\{\infty\}$. The \emph{strict $n$-subdivision} $S_n(X)$ is the category with objects $(i,x)$, where $x\in\Ob X$ and $0\le i\le n$ ($i\in\mathbb N$ if $n=\infty$), with
\[
\Hom\bigl((i,x),(j,y)\bigr)=
\begin{cases}
\Hom_X(x,y) & \text{if } i<j,\\
\{\id\} & \text{if } (i,x)=(j,y),\\
\emptyset & \text{otherwise},
\end{cases}
\]
and with the composition of $X$. It is acyclic, and it is finite if $X$ is finite and $n<\infty$. For a monoid $M$, regarded as a category with one object, we write $\CC_{M,n}=S_n(M)$: its objects are $0,\dots,n$, $\Hom(i,j)=M$ for $i<j$, and the composite of $s\co i\to j$ and $s'\co j\to l$ is $s's$. We write $p\co\CC_{M,n}\to M$ for the functor that sends each morphism to its label.
\end{definition}

Thus $\CC_{M,n}$ has exactly $n+1$ objects, and for $M=\Z/2$ it is the category $\CC_{\Z/2,n}$ of Subsection~\ref{subsec:RPn}.

\begin{lemma}\label{lem:Sn-sinks}
Let $X$ be a finite category and $n\ge1$. The objects $(n,x)$ are sinks of $S_n(X)$, there are no non-trivial morphisms between them, and the full subcategory on the remaining objects is $S_{n-1}(X)$. If they are added to $S_{n-1}(X)$ one at a time, in any order, the category of elements of $(n,x)$ is $S_{n-1}(X/x)$, where $X/x$ is the slice category, and its projection is $S_{n-1}$ of the forgetful functor $X/x\to X$.
\end{lemma}

\begin{proof}
The first assertions are immediate from the definition. When $(n,x)$ is added, no object $(n,y)$ maps to it, so the objects of its category of elements are the pairs $((i,y),u)$ with $i<n$ and $u\co y\to x$ in $X$, that is, the objects $(i,(y,u))$ of $S_{n-1}(X/x)$. A morphism $((i,y),u)\to((j,y'),u')$ is either an identity or, for $i<j$, a morphism $f\co y\to y'$ of $X$ with $u'\circ f=u$, that is, a morphism $(y,u)\to(y',u')$ of $X/x$.
\end{proof}

\begin{lemma}\label{lem:Sn-connected}
Let $Y$ be a finite category with a terminal object $\tau$, and let $n\ge0$. Then $BS_n(Y)$ is nonempty and $(n-1)$-connected.
\end{lemma}

\begin{proof}
We argue by induction on $n$. The category $S_0(Y)$ is discrete with object set $\Ob Y\neq\emptyset$. Let $n\ge1$. By Lemma~\ref{lem:Sn-sinks}, $S_n(Y)$ is obtained from $S_{n-1}(Y)$ by adding the sinks $(n,y)$, and we add $(n,\tau)$ first. Since $\tau$ is terminal, the forgetful functor $Y/\tau\to Y$ is an isomorphism. Hence so is the projection of the category of elements $S_{n-1}(Y/\tau)$ onto $S_{n-1}(Y)$, and by Lemma~\ref{lem:sink}(a) we obtain the cone on $BS_{n-1}(Y)$, which is contractible. Every other sink $(n,y)$ has category of elements $S_{n-1}(Y/y)$, whose projection lands in $S_{n-1}(Y)$, and $Y/y$ has the terminal object $(y,1_y)$, so $A_y=BS_{n-1}(Y/y)$ is nonempty and $(n-2)$-connected by induction. By Lemma~\ref{lem:sink}(b), $BS_n(Y)$ is obtained from the contractible space $\cone BS_{n-1}(Y)$ by attaching cones on the spaces $A_y$ along maps into it, so it is homotopy equivalent to $\bigvee_y\Sigma A_y$. The suspension of a nonempty $(n-2)$-connected space is $(n-1)$-connected.
\end{proof}

\begin{theorem}\label{thm:monoids}
Let $M$ be a finite monoid and $n\ge0$. Then $Bp\co B\CC_{M,n}\to BM$ is $n$-connected. In particular, $\CC_{M,n}$ is a finite acyclic category with $n+1$ objects that models $BM$ through degree $n$, $H_k(B\CC_{M,n};\Z)\cong H_k(BM;\Z)$ for $k<n$, and $\pi_1(B\CC_{M,n})\cong\pi_1(BM)$ if $n\ge2$.
\end{theorem}

\begin{proof}
Let $Y_M=M/{*}$ be the slice category of $M$ over its object. Its objects are the elements of $M$, its morphisms $s\co m\to m'$ are the elements $s\in M$ with $m's=m$, and its identity element $e$ is a terminal object, since $s\co m\to e$ forces $s=m$. By Lemma~\ref{lem:Sn-sinks}, $\CC_{M,k+1}$ is obtained from $\CC_{M,k}$ by adding the sink $k+1$, whose category of elements is $S_k(Y_M)$. By Lemma~\ref{lem:sink}(b) and Lemma~\ref{lem:Sn-connected}, $B\CC_{M,k+1}$ is obtained from $B\CC_{M,k}$ by attaching the cone on a nonempty $(k-1)$-connected space. Replacing that space by a homotopy equivalent CW complex with one $0$-cell and no other cells of dimension less than $k$, we see that the pair $(B\CC_{M,k+1},B\CC_{M,k})$ is $k$-connected \cite[Section~4.1]{Hatcher2002}.

The nerve of $\CC_{M,\infty}$ is the union of the nerves of the $\CC_{M,k}$. Since spheres and their homotopies are compact, the inclusion $B\CC_{M,n}\to B\CC_{M,\infty}$ is $n$-connected. It remains to see that $p\co\CC_{M,\infty}\to M$ induces a homotopy equivalence. The comma category $p/{*}$ has objects $(i,m)$ with $m\in M$, and its morphisms $(i,m)\to(j,m')$, for $i<j$, are the $s\in M$ with $m's=m$. Hence $p/{*}\cong S_\infty(Y_M)$. Its classifying space is the union of the spaces $BS_k(Y_M)$, which are $(k-1)$-connected by Lemma~\ref{lem:Sn-connected}. So it is weakly contractible, hence contractible, and Quillen's Theorem~A \cite{Quillen1973} shows that $Bp$ is a homotopy equivalence. The statements about homology and fundamental groups follow from the properties of $n$-connected maps recalled in Subsection~\ref{subsec:nu}.
\end{proof}

\subsection{Finite groups}\label{subsec:groups}

For groups, the categories $\CC_{G,n}$ are quotients of posets, as in Section~\ref{sec:projective}, and their number of objects cannot be improved.

\begin{corollary}\label{cor:groups}
Let $G$ be a finite group and $n\ge0$.
\begin{enumerate}[(a)]
\item $G$ acts freely on the objects of the stack $\PP_{G,n}$ of $n+1$ copies of the discrete category with set of objects $G$, by $g\cdot(i,a)=(i,ag^{-1})$, and $\PP_{G,n}/G\cong\CC_{G,n}$. Consequently $B\CC_{G,n}$ is homeomorphic to the quotient of the join $G^{*(n+1)}$ by the diagonal action of $G$.
\item The map $B\CC_{G,n}\to BG$ is $n$-connected. In particular $\pi_1(B\CC_{G,n})\cong G$ for $n\ge2$.
\end{enumerate}
\end{corollary}

\begin{proof}
(a) The orbit of the morphism $(i,a)\to(j,b)$ is determined by $i<j$ and by $ba^{-1}\in G$, and if $(i,a)\to(j,b)$ and $(j,b)\to(l,c)$ are composable, then $(cb^{-1})(ba^{-1})=ca^{-1}$. By Lemma~\ref{lem:free-quotient}, $\PP_{G,n}/G\cong\CC_{G,n}$ and $B\CC_{G,n}\cong B\PP_{G,n}/G$. By Lemma~\ref{lem:join}(a), $B\PP_{G,n}$ is the join of $n+1$ copies of $G$, and the homeomorphism is equivariant, since it is defined simplex by simplex. (b) This is Theorem~\ref{thm:monoids}.
\end{proof}

For $G=\Z/2$ we recover Theorem~\ref{thm:RPn} and the inclusion of $\RP^n$ in $\RP^\infty$.

\begin{example}\label{ex:perfect}
For $n=2$, Corollary~\ref{cor:groups} shows that three objects already produce every finite group as a fundamental group, and the computation can be made by hand. Write the non-trivial morphisms of $\CC_{G,2}$ as $f_u\co0\to1$, $g_v\co1\to2$ and $h_w\co0\to2$, so that $g_v\circ f_u=h_{vu}$. With the maximal tree $\{f_e,g_e\}$, the relations $[f_u][g_v]=[h_{vu}]$ give $[g_v]=[h_v]$ and $[f_u]=[h_u]$, and the classes $z_u=[h_u]$ satisfy $z_e=1$ and $z_uz_v=z_{vu}$. This is a presentation of the opposite group of $G$, which is isomorphic to $G$ by $u\mapsto u^{-1}$. Since $B\CC_{G,2}$ is $2$-dimensional with Euler characteristic $3-3|G|+|G|^2$, and $H_2$ is free by Lemma~\ref{lem:faces}(c), we get $H_1\cong G^{\mathrm{ab}}$ and $H_2\cong\Z^{(|G|-1)(|G|-2)}$. For the alternating group $A_5$, which is perfect, $B\CC_{A_5,2}$ has the homology of a bouquet of $3422$ spheres of dimension $2$, but its fundamental group is $A_5$. Thus homology no longer determines the homotopy type of a three-object category when there are more than two morphisms between the first two objects: here the target $B\K_{|G|}$ of $p_T$ in \eqref{eq:three} is a wedge of circles, and Proposition~\ref{prop:three} does not apply.
\end{example}

\begin{proposition}\label{prop:groups-minimal}
Let $G$ be a nontrivial finite group. Then $\nu_n(BG)=n+1$ for every $n\ge0$, and $\CC_{G,n}$ realizes it.
\end{proposition}

\begin{proof}
The upper bound is Corollary~\ref{cor:groups}(b), and the case $n=0$ is clear. Let $n\ge1$ and let $\CC$ be a finite acyclic category with an $n$-connected map $f\co B\CC\to BG$; we show that $|\Ob\CC|\ge n+1$. If $n=1$, then $\pi_1(B\CC)$ maps onto $G\neq1$, so $\CC$ has at least two objects. Let $n\ge2$. Choose a prime $\ell$ dividing $|G|$ and a subgroup $H\le G$ of order $\ell$. Let $Z\to BG$ be the covering with fibre $G/H$, so that $Z\simeq BH$, and let $\widehat Y\to B\CC$ be its pullback along $f$. Since $f$ is an isomorphism on $\pi_1$, $\widehat Y$ is connected. Comparing the long exact sequences of the two coverings shows that the induced map $\widehat Y\to Z$ is again $n$-connected. Since $\widehat Y$ is a finite covering of the trisp $B\CC$, it is a trisp of the same dimension $d\le|\Ob\CC|-1$ (Lemma~\ref{lem:faces}(a)), and its top homology group $H_d(\widehat Y;\Z)$, being a group of cycles in the free group of $d$-chains, is free. Recall that $H_k(BH;\Z)\cong\Z/\ell$ for odd $k$ \cite{Hatcher2002}. If $n$ is odd, $H_n(\widehat Y;\Z)\to H_n(BH;\Z)\cong\Z/\ell$ is onto, so $d\ge n$. If $n$ is even, $H_{n-1}(\widehat Y;\Z)\cong\Z/\ell$ has torsion, so $d>n-1$. In both cases $|\Ob\CC|\ge d+1\ge n+1$.
\end{proof}

Since $\nu(BG)\ge\nu_n(BG)=n+1$ for every $n$, this recovers $\nu(BG)=\infty$, now for every nontrivial finite group, and shows that the sequence $\nu_n$ need not stabilize. Real projective spaces show both behaviours at once.

\begin{example}\label{ex:nu-n-RP}
Let $k\ge1$. Then $\nu_n(\RP^k)=\min\{n,k\}+1$. Indeed, let $n\le k$. The inclusion $\RP^n\subseteq\RP^k$ is $n$-connected, so $\CC_{\Z/2,n}$ gives $\nu_n(\RP^k)\le n+1$ by Theorem~\ref{thm:RPn}. Conversely, the inclusion $\RP^k\subseteq\RP^\infty=B\Z/2$ is $k$-connected, so composing it with an $n$-connected map $B\CC\to\RP^k$ gives an $n$-connected map to $B\Z/2$, and $|\Ob\CC|\ge n+1$ by Proposition~\ref{prop:groups-minimal}. For $n\ge k$ we get $k+1=\nu_k(\RP^k)\le\nu_n(\RP^k)\le\nu(\RP^k)=k+1$. Thus $\nu_n(\RP^k)$ follows $\nu_n(\RP^\infty)=n+1$ up to $n=k$ and then stays at $\nu(\RP^k)$.
\end{example}

The lens models of Section~\ref{sec:projective} fit into the same picture: sending $A_i$ to $2i$ and $B_i$ to $2i+1$ identifies $\mathcal L_m(q_0,\dots,q_k)$ with the wide subcategory of $\CC_{\Z/m,2k+1}$ that keeps, between $2i$ and $2i+1$, only the morphisms labeled $0$ and $\bar q_i$. Both categories have $2k+2$ objects and model $B\Z/m$ through degree $2k+1$, the lens models with fewer morphisms and with a manifold as classifying space.

\subsection{Finite simply connected complexes}\label{subsec:approx}

The following consequence depends on a recent result of Ben-Bassat, who proved Fiedorowicz's conjecture that every finite simply connected CW complex is homotopy equivalent to the classifying space of a finite monoid \cite[Theorem~C]{BenBassat2026}.

\begin{corollary}[{via \cite[Theorem~C]{BenBassat2026}}]\label{cor:approx}
Let $X$ be a finite simply connected CW complex. Then $\nu_n(X)\le n+1$ for every $n\ge0$: there is a finite acyclic category $\CC$ with $n+1$ objects and an $n$-connected map $B\CC\to X$, so that $H_k(B\CC;\Z)\cong H_k(X;\Z)$ and $\pi_k(B\CC)\cong\pi_k(X)$ for $k<n$.
\end{corollary}

\begin{proof}
By \cite[Theorem~C]{BenBassat2026}, $X\simeq BM$ for a finite monoid $M$. Take $\CC=\CC_{M,n}$ and apply Theorem~\ref{thm:monoids}.
\end{proof}

\begin{figure}[ht]
\centering
\begin{tikzpicture}[>={Stealth[length=4pt]}, lbl/.style={font=\scriptsize}]
\node (0) at (0,0) {$0$};
\node (1) at (2.3,2.0) {$1$};
\node (2) at (4.6,0) {$2$};
\foreach \b in {-22,-11,11,22} {
  \draw[->, gray!70] (0) to[bend left=\b] (1);
  \draw[->, gray!70] (1) to[bend left=\b] (2);
  \draw[->, gray!70] (0) to[bend left=\b] (2);
}
\draw[->, thick] (0) -- node[lbl, above left]{$x_{12}$} (1);
\draw[->, thick] (1) -- node[lbl, above right]{$x_{21}$} (2);
\draw[->, thick] (0) -- node[lbl, below=2pt]{$x_{22}$} (2);
\end{tikzpicture}
 
\smallskip
\caption{The category $\CC_{M,2}$, with five morphisms between any two objects, labeled by $M$. The thick morphisms form a commutative triangle, since the composite of $x_{12}\co0\to1$ and $x_{21}\co1\to2$ is $x_{21}x_{12}=x_{22}$.}
\label{fig:band}
\end{figure}
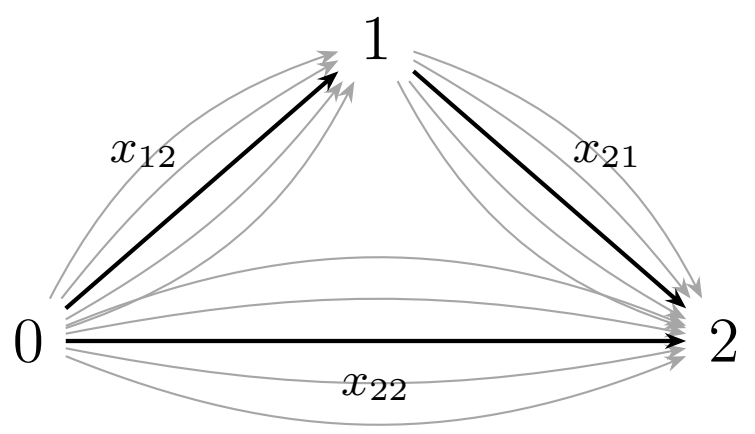

\begin{example}\label{ex:band}
Let $M=\{1,x_{11},x_{12},x_{21},x_{22}\}$ be the $2\times2$ rectangular band with an identity adjoined, with product $x_{ij}x_{kl}=x_{il}$. Fiedorowicz showed that $BM\simeq S^2$ \cite{Fiedorowicz2002}. For $n=2$, the category $\CC_{M,2}$ has three objects and five morphisms between any two of them, labeled by $M$; see Figure~\ref{fig:band}. By Theorem~\ref{thm:monoids}, $B\CC_{M,2}$ is simply connected and maps onto $\pi_2(S^2)$, but it is not yet a sphere: it is a bouquet of twelve $2$-spheres. As $n$ grows, $B\CC_{M,n}$ agrees with $S^2$ in more and more degrees, and in the limit $B\CC_{M,\infty}\simeq S^2$ by Corollary~\ref{cor:paradox}. The sphere category $\Sph$ does better with the same three objects, since $\nu_n(S^2)=3$ for $n\ge2$ (Remark~\ref{ex:nu-n-spheres}): an exact model is a matter of choosing the right composition table, not of adding objects.
\end{example}

\begin{remark}\label{rem:monoids-general}
For a finite simply connected CW complex $X$, Corollary~\ref{cor:approx} and Proposition~\ref{prop:stabilization} show that $\nu_n(X)\le n+1$ for every $n$ and that $\nu_n(X)=\nu(X)$ for $n$ large. Below a given degree, the complexity of $X$ can always be stored in the morphisms, through the size of $M$; what $\nu(X)$ measures is the cost of killing, with finitely many objects, the homotopy that $B\CC_{M,n}$ has in degrees $\ge n$ and $X$ does not. The previous example illustrates this point.
\end{remark}

Letting $n$ go to infinity gives a statement that sounds paradoxical.

\begin{corollary}\label{cor:paradox}
Every connected CW complex $X$ is homotopy equivalent to $B\CC$ for an acyclic category $\CC$ with countably many objects. If $X$ is finite and simply connected, $\CC$ can be chosen with finite hom-sets.
\end{corollary}

\begin{proof}
By McDuff's theorem \cite{McDuff1979}, $X\simeq BM$ for a discrete monoid $M$, and by \cite[Theorem~C]{BenBassat2026} $M$ can be taken finite when $X$ is finite and simply connected. Take $\CC=\CC_{M,\infty}$. As in the proof of Theorem~\ref{thm:monoids}, $p/{*}\cong S_\infty(Y)$ with $Y=M/{*}$, and $Y$ has the terminal object $e$. This category is contractible for any monoid, finite or not. Let $D(i,y)=(2i+1,y)$ and $T(i,y)=(2i+2,e)$, with $D(f)=f$ and $T(f)=\id_e$. There are natural transformations $1\Rightarrow D$, with components $\id_y$; $D\Rightarrow T$, with components the unique morphisms $y\to e$; and $\mathrm{const}_{(0,e)}\Rightarrow T$, with components $\id_e$. So the identity is homotopic to a constant functor, and by Quillen's Theorem~A \cite{Quillen1973}, $Bp\co B\CC_{M,\infty}\to BM$ is a homotopy equivalence.
\end{proof}

The paradox dissolves once one notices that $\CC_{M,\infty}$ carries no more information than $M$: the infinitely many objects only unfold the endomorphisms of a single object, and the excess homotopy that $B\CC_{M,n}$ has in degree $n$ is killed at the next stage.

\section{Comparison and final remarks}\label{sec:comparison}

We first compare our models with minimal finite models. A finite model of $S^n$ needs $2n+2$ points \cite[Theorem~2.13]{BarmakMinian2007}, against $\nu(S^n)=n+1$. A minimal finite model of $\bigvee^rS^1$ has $\min\{i+j:(i-1)(j-1)\ge r\}$ points \cite{BarmakMinian2007}, a number that grows like $2\sqrt r$, against $\nu=2$. A finite model of $\RP^2$ needs $13$ points, and there are exactly two with $13$ points \cite[Corollary~5.8, Theorem~5.11]{CianciOttina2018}, against $\nu(\RP^2)=3$; an open cover of $\RP^2$ gives an acyclic model with seven objects \cite{CarcaciaCampos2026OpenCovers}. More generally, torsion in the homology of a finite poset forces at least $13$ points \cite[Theorem~5.10]{CianciOttina2018}, and the finite models of $\RP^n$ and of $\M(\Z/m,1)$ constructed by Cianci and Ottina have $(3^{n+1}-1)/2$ and $4m+5$ points \cite[Propositions~7.1 and~7.7]{CianciOttina2018}, against $n+1$ and $3$ objects. For closed PL $n$-manifolds, a crystallization encodes a pseudo-triangulation with exactly $n+1$ vertices \cite{Pezzana1974,Pezzana1975}, as many as the objects of our models of spheres, real projective spaces and lens spaces; the difference is that our models are nerves of categories, in which every chain of composable morphisms spans a simplex.

Behind all these computations there is a single mechanism. Choosing a linear extension of $\preceq$, every finite acyclic category is obtained by adding one sink at a time, and by the sink lemma each step attaches a cone along an explicit map of nerves: the sink lemma is the categorical form of attaching a cell. The analogy is not exact, and the difference is what $\nu$ measures. A cell can be attached along any map, whereas a sink $t$ is attached along $B(\DD/t)$, and $\dim\nerve(\DD/t)\le\dim\nerve(\DD)$, since the projection sends non-trivial morphisms to non-trivial morphisms; thus each object raises the dimension by at most one. For real projective spaces, lens spaces and the Moore spaces $\M(\Z/m,1)$, whose minimal cell structures have one cell in each dimension, this costs nothing, and our models have one object for each cell. For spheres it costs $n-1$ extra objects, which the parallel morphisms then exploit to carry whole bouquets at no further cost in objects.

The idea that runs through this paper is one of the oldest in algebraic topology: to understand a space by building it from simple pieces glued along explicit maps. It goes back to Poincar\'e, who studied manifolds through their decompositions into cells and simplices \cite{Poincare1895,Poincare1899}, and it lives on in simplicial complexes, simplicial sets and CW complexes. Finite acyclic categories are one more incarnation of it, in which each new object is the apex of a cone and the morphisms into it are the instructions for gluing it. What changes from one incarnation to another is where the information is stored. In a finite space two points are related at most once, and complexity has to be paid for in points. In a monoid it is concentrated in the endomorphisms of a single object \cite{McDuff1979,BenBassat2026}. Acyclic categories lie in between, and store it in parallel morphisms and their composition: odd spheres in joins of Kronecker categories, the Moore spaces $\M(\Z/m,1)$ in a table that adds modulo $m$, real projective spaces and lens spaces in morphisms labeled by a cyclic group. The two invariants measure how far this middle ground reaches. Below a fixed degree it reaches as far as monoids do, since spreading a monoid over $n+1$ objects keeps its homotopy below degree $n$; matching a space exactly requires finitely many objects to absorb all of its homotopy at once. The lens spaces $L(m;1,\dots,1)$ of increasing dimension illustrate the difference: each has an object-minimal model, and their union $S^\infty/(\Z/m)$ is a model of $B\Z/m$, which has none. Classifying spaces of finite groups show at the same time the reach and the limit of finite acyclic models: for every nontrivial finite group $G$,
\[
\nu_n(BG)=n+1\quad\text{for every }n\ge0,\qquad\nu(BG)=\infty.
\]
Each new degree of $BG$ costs exactly one more object, and no finite number of objects sees all of them.

\bibliographystyle{amsplain}
\bibliography{minimal_acyclic_models_clean}

\end{document}